%% file: main1.tex
\documentclass[a4paper,10pt]{article}
\usepackage[utf8]{inputenc}

\providecommand{\keywords}[1]{\textbf{\textit{keywords}} #1}

\usepackage{amsthm}
\usepackage{amsmath}
\usepackage{amsfonts}
\usepackage{amssymb}
\usepackage{epsfig,epstopdf,titling}
\usepackage{caption}
\usepackage{subcaption}
\usepackage{multirow}
\usepackage{url}
\usepackage{graphicx}
\usepackage{pgfplots}
\usepackage{array,tabularx}
\usepackage{multicol}  
\usepackage{algorithmic} 
\usepackage[boxruled,longend,linesnumbered,ruled]{algorithm2e}
\usepackage{subcaption}
\theoremstyle{definition}

\newtheorem{example}{Example}[section]

\theoremstyle{definition}
\usepackage{mathrsfs}

\theoremstyle{remark}

\usepackage{authblk}

\title{Frequency Limited $\mathcal{H}_2$ Optimal Model Reduction of Large-Scale Sparse Dynamical Systems}

\author[]{Xin Du\thanks{ School of Mechatronic Engineering and Automation, Shanghai University, Shanghai-200072, China
 and   
Key Laboratory of Modern Power System Simulation and Control \& Renewable Energy Technology, Ministry of Education(Northeast Electric Power University), Jilin-132012, China
, {duxin@shu.edu.cn}}}
\author[]{M. Monir Uddin\thanks{Department of Mathematics and Physics, North south University, Dhaka-1229, Bangladesh, {monir.uddin@northsouth.edu}}}
\author[]{A. Mostakim Fony\thanks{Department of Mathematics, Chittagong University,  Chittagong, Bangladesh, {asibmostakim1995@gmail.com}}}
\author[]{Md. Tanzim Hossain\thanks{Department of Electrical and Computer Engineering, North South University, Dhaka-1229, Bangladesh, {tanzim.hossain@northsouth.edu}}}
\author[]{Mohammaed Sahadat-Hossain\thanks{Department of Mathematics and Physics, North south University, Dhaka-1229, Bangladesh, {mohammad.hossain@northsouth.edu}}}
\affil[]{}
\date{}
\begin{document}
  \maketitle
\begin{abstract}
We mainly consider the frequency limited $\mathcal{H}_2$ optimal model order reduction of large-scale sparse generalized systems. For this purpose we need to solve two Sylvester equations. This paper proposes efficient algorithm to solve them efficiently. The ideas are also generalized to index-1 descriptor systems. Numerical experiments are carried out using Python Programming Language and the results are presented to demonstrate the approximation accuracy and computational efficiency of the proposed techniques. 
\end{abstract}
\begin{center}
\end{center}
\section{Introduction}\label{sec:introduction}
Model order reduction (MOR) is a process to approximate a high-order dynamical system by a substantially low-order system with a maximum accuracy. This tool is now widely used in different disciplines of science, engineering and technology to reduce the complexity of the model. In general, the reduced order models are used in controller design, simulation and optimization. For motivation, applications and techniques of MOR see, e.g,. \cite{morAnt05,morUdd19}.

The commonly useful methods for the model reduction of large-scale linear time invariant dynamical systems are the  balanced truncation and $\mathcal{H}_2$ optimal model reduction \cite{morAnt05}.  Both the methods are well established and  successfully investigated to find the model reduction of large-scale sparse  dynamical systems.  In the recent times frequency and time limited model reduction methods have taken  a lot of attentions due to its demand in real-life applications. In many applications, a specific frequency interval is more important, i.e., the ROM should maintain a superior accuracy within that desired frequency interval. Balanced truncation based frequency limited model reduction was discussed by Gawronski and Juang in \cite{morGawJ90}.  The computational techniques for the time and frequency limited balanced truncation were discussed in \cite{morBenKS16}. 

The optimal $\mathcal{H}_2$ model reduction methods have been studied and investigated  in  \cite{morVanGA08,morGugAB08,morHylB85,morYanJ99,morXuZ11}. See, also the references cited therein. In all these papers the proposed technique is based on either Gramian assistance first-order optimality conditions \cite{morVanGA08} or the tangential interpolaton \cite{morGugAB08} of the transfer function. In fact both the conditions are coincided which is shown in \cite{morXuZ11}. These papers only discuss the model reduction on the infinite frequency interval.  For the time limited case the we refer the readers to \cite{morgoyR19}. Although, in \cite{morpetL14} authors briefly introduced optimal $\mathcal{H}_2$ model reduction problem of standard state space systems considering a restricted frequency interval, there the implementation details were not given.

This paper focuses on the computational techniques of the frequency limited optimal $\mathcal{H}_2$ model reduction method of large-scale sparse systems. We mainly generalized the idea as in \cite{van2008h2,morXuZ11} in which the proposed algorithm is called two sided iteration algorithm (TSIA). Moreover, to implement the frequency limited TSIA we need to solve two frequency limited Sylvester equations. This paper also discusses how to solve the Sylvester equations efficiently by preserving the sparsity of the system. Besides the generalized systems the idea is also extended for index-1 descriptor systems. The benefits of the algorithmic improvements presented in this paper are illustrated by several numerical examples. We have generated the the results by using Python Programming Language. Rest of this paper is organized as follows.

Section~\ref{sec:gensys} overview  the TSIA and the optimal $\mathcal{H}_2$ model reduction of generalized system. Then the ideas of this section are discussed in the next sections  for  the frequency-limited model order reduction. Section~\ref{sec:flmor-1:index1sys} presents the algorithm for solving frequency limited Sylvesters equations which provides projectors to carry out the FLMOR. The results of the numerical experiments  are depicted in Section~\ref{sec:numreslt} which show the efficiency and capability of the proposed methods. 

\section{The TSIA and $\mathcal{H}_2$ optimal model order reduction of generalized Systems}
\label{sec:gensys}
The goal of this section is to review the basic idea of $\mathcal{H}_2$ optimal model order reduction of generalized Systems. Let us consider a linear time invariant continuous-time system of the form
\begin{equation}\label{eq:standard:ltisystem}
\begin{aligned}
E\dot x(t) & = A x(t)+B u(t),\\
      y(t) & = C x(t)+ D u(t),
\end{aligned}
\end{equation}
where $ E \in \mathbb{R}^{n\times n}$ is non-singular, 
and $ A\in\mathbb R ^{n\times n}$, $B\in\mathbb R^{n\times p}$,
$C\in\mathbb R ^{m\times n}$ and ${D}\in\mathbb R^{m\times p}$. 
The transfer-function matrix of this system is defined by 
$\text{G}(s)= C(sE-A)^{-1}B + D$, where $s\in \mathbb{C}$. 
The controllability and the observability Gramians of the system on the infinite frequency range can be defined as \cite{morZhoSW99}
\begin{align*}
  P & = \frac{1}{2\pi}\int_{-\infty}^{\infty}(\imath \omega E
  -A)^{-1}B B^{T}(\imath \omega E^{T} - A^{T}) \mathop{d}\omega \quad \text{and}\\
  Q & = \frac{1}{2\pi}\int_{-\infty}^{\infty}(\imath \omega E^{T}
  -A^{T})^{-1}C^{T}C(\imath \omega E - A) \mathop{d}\omega,
\end{align*}
and they are the solutions of the continuous-time algebraic Lyapunov equations
\begin{equation}
  \begin{aligned}
    \label{eq:stablyapunov}
    & APE^T+EPA^T+BB^T = 0 \quad \text{and} \\ 
    & A^TQE+E^TQA+C^TC = 0,  
  \end{aligned}
\end{equation}
respectively. We want to construct a substantially reduced-order model
\begin{equation}
\label{eq:standard:reducedltisystem}
\begin{aligned}
\dot {\hat{x}}(t) & = \hat{A} \hat{x}(t)+\hat{B} u(t),\\ 
              \hat{y}(t) & = \hat{C}\hat{x}(t)+ \hat{D} u(t),
\end{aligned}
\end{equation}
where $ \hat{A}\in\mathbb R ^{r\times r}$, $\hat{B}\in\mathbb R^{r\times p}$, $\hat{C}\in\mathbb R ^{m\times r}$ and $\hat{D}\in\mathbb R^{m\times p}$. 
The goal is to minimize the error 
\begin{align}\label{eq:gensys:error}
\Xi = \| \text{G}-\hat{\text{G}}\|_{\mathcal{H}_2},
\end{align}
where $\| . \|_{\mathcal{H}_2}$ denotes the system's $\mathcal{H}_2$-norm \cite{morAnt05} and $\hat{G}(s) = \hat{C}(sI-\hat{A})^{-1}\hat{B} + \hat{D}$ is 
the transfer-function matrix of the reduced  system.  
The $\mathcal{H}_2$-norm of the error system as defined in (\ref{eq:gensys:error}) can be measured  by 
\begin{align}\label{eq:gensys:eerr}
 \Xi = : \sqrt{\mathrm{Tr}\left(C_{\Xi}P_{\Xi}C_{\Xi}^T\right)} \quad \text{or}\quad \sqrt{\mathrm{Tr}\left(B_{\Xi}^TQ_{\Xi}B_{\Xi}\right)} ,
\end{align}
where $P_\Xi$ and $Q_\Xi$ are the solutions of the Lyapunov equations %
\begin{subequations} \label{eq:gensys:rsleq}
  \begin{align}
   & A_\Xi P_\Xi E_\Xi^{T} + E_\Xi P_\Xi A_\Xi^{T}  + B_\Xi B_\Xi^{T} =0  \quad \text{and} \label{eq:gensys:rsleq1} \\
   & A_\Xi^{T} Q_\Xi E_\Xi + E_\Xi^{T} Q_\Xi A_\Xi  + C_\Xi^{T}C_\Xi=0,\label{eq:gensys:rsleq2}
  \end{align}
\end{subequations}
in which 
\begin{align*}
E_\Xi = \begin{bmatrix} E & \\ & \hat{I}\end{bmatrix}, \quad
A_\Xi = \begin{bmatrix} A & \\ & \hat{A}\end{bmatrix}, \quad
B_\Xi = \begin{bmatrix} B  \\  \hat{B}\end{bmatrix} \quad \text{and}\quad 
C_\Xi = \begin{bmatrix} C  & -\hat{C}\end{bmatrix}.
\end{align*}
Now, partitioning $P_\Xi$ and $Q_\Xi$ as 
\begin{align*}
 P_\Xi = \begin{bmatrix} P & M \\ M^T & \hat{P}\end{bmatrix} \quad
\text{and} \quad Q_\Xi = \begin{bmatrix} Q & N \\ N^T & \hat{Q}\end{bmatrix},\, \text {respectively} 
\end{align*}
and plugging into (\ref{eq:gensys:rsleq}) we obtain the following algebraic matrix equations
\begin{subequations}
  \label{eq:gensys:splitedlyap}
  \begin{align}
 A P E^{T} + E P A^{T} + B B^{T} & = 0,  \label{eq:gensys:splitedlyap1} \\
 \hat{A} \hat{P}  +  \hat{P} \hat{A}^{T}  +\hat{B} \hat{B}^{T} & = 0, \label{eq:gensys:splitedlyap2} \\
A M  + E M \hat{A}^{T} + B \hat{B}^{T} & =0,  \label{eq:gensys:splitedlyap3} \\
 A^{T} Q E + E^{T} Q A  + C^{T}C & =0,\label{eq:gensys:splitedlyap4} \\
\hat{A}^{T} \hat{Q}  +  \hat{Q} \hat{A}  + \hat{C}^{T}\hat{C} & = 0,\label{eq:gensys:splitedlyap5} \\
A^{T} N + E^{T} N \hat{A}  - C^{T}\hat{C} & = 0,\label{eq:gensys:splitedlyap6} 
\end{align}
\end{subequations}
where $\hat{P}$ and $\hat{Q}$ are respectively known as the controllability and observability Gramians of the reduced systems. Therefore, the $\mathcal{H}_2$ norm of the error system in (\ref{eq:gensys:eerr}) can be measured by 
\begin{equation}\label{eq:gensys:eerrex}
 \Xi = :  \begin{cases}\sqrt{{\mathrm{Tr}\left(CPC^T\right)} +{\mathrm{Tr}\left(\hat{C}\hat{P}\hat{C}^T\right)}+2{\mathrm{Tr}\left(CM\hat{C}^T\right)}}\\ 
 \quad \text{or}\quad \\
 \sqrt{{\mathrm{Tr}\left(B^TQB\right)} +{\mathrm{Tr}\left(\hat{B}^T\hat{Q}\hat{B}\right)}+2{\mathrm{Tr}\left(B^TN\hat{B}\right)}}.
          \end{cases}
\end{equation}
The first-order optimality conditions for the optimal $\mathcal{H}_2$ model reduction was given in \cite{morWil70}, which is known as  Wilson conditions. Wilson conditions are based on the first derivatives of (\ref{eq:gensys:error}) with respect to $\hat{A}$, $\hat{B}$ and $\hat{C}$ as follows:
\begin{align*}
 \nabla\Xi_{\hat{A}}= 2(\hat{Q}\hat{P}+W^TEV),\, 
 \nabla\Xi_{\hat{B}}= 2(\hat{Q}\hat{B}+W^TE^{-1}B),  \,
 \nabla\Xi_{\hat{C}}= 2(\hat{C}\hat{P}-CV).
\end{align*}
Setting these three derivatives to zero leads to 
the \emph{Wilson conditions},
\begin{align}
 \hat{Q}\hat{P}+N^T E M = 0,\label{eq:gensys:wilsoncondition1}\\
 \hat{Q}\hat{B}+N^TE^{-1}B=0,\label{eq:gensys:wilsoncondition2} \\
 \hat{C}\hat{P}-CM=0.\label{eq:gensys:wilsoncondition3} 
\end{align}
These three conditions in fact yield the left and and right projection matrices to compute an optimal reduced order system (\ref{eq:standard:reducedltisystem}) 
and in the  optimal reduced order system the reduced matrices are formed
as 
\begin{align}\label{eq:gensys:reducedmatrices}
    \hat{A}= W^TE^{-1}A V,\quad \hat{B}= W^T B,\quad\,\text{and}\quad \hat{C}= C V,
\end{align}
where $V=M\hat{P}^{-1}$ and $W^T=-\hat{Q}^{-1}N^T$ and hence it can be proved that $W^TEV=I$. However, we can not guarantee that $\hat{P}$ and $\hat{Q}$ are invertible, since to assure this the reduced model should be completely controllable and observable \cite{morAnt05}. In the case that they are invertible, the multiplication from the right is only a transformation of bases and does not change the subspace. The idea by Xu and Zeng \cite{morXuZ11} was to satisfy the Wilson conditions by setting 
\begin{align}
W = N \quad\text{and}\quad V = M.
\end{align}
Note that $V$ and $W$ can be computed by solving the matrix equations (\ref{eq:gensys:splitedlyap3}) and (\ref{eq:gensys:splitedlyap6}), respectively which are known as Sylvesters equations. Another important observation is that if we want to compute the optimal projection subspace we already need the optimal solution $\hat{A}$, $\hat{B}$ and $\hat{C}$. However, this is not known prior. A possible solution is to start with a reduced model, which emerged from an arbitrary projection of the original model, solve matrix equations (\ref{eq:gensys:splitedlyap3}) and (\ref{eq:gensys:splitedlyap3}), compute the projectors, and restart the process with the newly obtained reduced model until we are satisfied. In this way we get a kind of a fixed point iteration. This procedure is called two sided iteration algorithm (TSIA) by Xu and Zeng in \cite{morXuZ11}. 


The Wilson conditions are Gramian-based conditions since they are related to Gramians of the system. The Hyland-Bernstein conditions \cite{morHylB85} are another gramian-based first-order optimal conditions, which were shown to be equivalent to the Wilson conditions \cite{gugercin2008h_2}. Van Dooren et al., were characterized the tangential interpolation based  $\mathcal{H}_2$ optimal conditions in \cite{van2008h2}. One drawback of interpolation based  model reduction is to selection of interpolation points. However in \cite{gugercin2008h_2} authors proposed the Iterative Rational Krylov Algorithms (IRKA) to resolve this problem. On the other hand Xu and Zeng  showed in \cite{morXuZ11} that both the Gramian and interpolation based optimality conditions are the same. In \cite{morXuZ11} authors also presented the two sided iterative algorithm (TSIA) for a standard system which is slightly modified as in \cite{van2008h2}. For our convenient the TSIA for the $\mathcal{H}_2$ optimal model reduction for the generalized system (\ref{eq:standard:ltisystem}) is summarized in Algorithm \ref{alg:irka1} 
\begin{algorithm}[t]
	\SetAlgoLined
	\SetKwInOut{Input}{Input}
	\SetKwInOut{Output}{Output}
	\caption{Two-sided iteration algorithm (TSIA).}
	\label{alg:irka1}
	\Input {$E, A, B, C, D$.}
	\Output{$\hat{A}, \hat{B}, \hat{C}$, $\hat{D}: = D$.}
	Choose matrices $W_0\in \mathbb{R}^{n\times r}$ and $V_0\in\mathbb{R}^{n\times r}$  such that $W_0^TV_0=I$.\\
	Construct the reduced-order matrices:
    $$\hat{A} = W_0^T E^{-1} A V_0,\, \hat{B} = W_0^TE^{-1}B\quad \text{and}\quad \hat{C} = C V_0$$.\\
	\While{$i\leq N-1$}{%
	Compute $V_i$ and $W_i$ by solving Sylvesters 
	\begin{subequations}
  \label{eq:gensys:splitedlyapex}
  \begin{align}
& A V  + E V \hat{A}^{T} + B \hat{B}^{T}=0  \label{eq:gensys:splitedlyapex3e} \\
&  A^{T} W + E^{T} W \hat{A}  - C^{T}\hat{C}=0,\label{eq:gensys:splitedlyap6e} 
\end{align}
\end{subequations}

	Compute $W_{i+1}=W_i(V_i^TW_i)^{-1}$ and $V_{i+1}=V_i$ \\
	Construct the reduced-order matrices
	$\hat{A} = W_{i+1}^T E^{-1}A V_{i+1}, \hat{B} = W_{i+1}^TE^{-1}B$ and $\hat{C} = C V_{i+1}$.\\
	$i=i+1$. \\}
	%
\end{algorithm}
%
\section{FL $\mathcal{H}_2$ optimal MOR of generalized systems}\label{sec:flmor}
This section moves to the frequency limited $\mathcal{H}_2$ optimal model reduction of the system (\ref{eq:standard:ltisystem}). For this purpose we first define frequency limited Gramians. In the previous section the system Gramians have been  defined on the infinite frequency interval. If we replace  the infinite interval (-$\infty$, $\infty$) into a finite interval $\omega=[\omega_1$, $\omega_2]$) then the controlability and the observablity Gramians can be defined as
\begin{align}\label{eq:sec:flmor:flgram}
  P_\omega &=\frac{1}{2\pi}\int_{\omega_1}^{\omega_2}(\imath \nu E -A)^{-1}B B^{T}(\imath \nu E^{T} - A^{T}) \mathop{d}\nu \\
  Q_\omega &=\frac{1}{2\pi}\int_{\omega_1}^{\omega_2}(\imath \nu E^{T} -A^{T})^{-1}C^{T}C(\imath \nu E - A) \mathop{d}\nu,
\end{align}
which satisfies the frequency limited Lyapunov equations

 \begin{equation}\label{eq:flmor:flimgrams}
  \begin{aligned}
    AP_\omega E^T+E P_\omega A^T+B_\omega BB^T+BB^TB_\omega^\ast=0,  \\ 
    A^T Q_\omega E+E^T Q_\omega A+C_\omega^\ast C^TC+C^TC C_\omega=0,  
  \end{aligned}
\end{equation}
with 
\begin{align*}
 & B_\omega = \frac{\imath}{2\pi}\ln(A+\imath\omega_2 E)(A+\imath\omega_1 E)^{-1}\,\,\text{and}\\
 & C_\omega = \frac{\imath}{2\pi}\ln(A+\imath\omega_1 E)^{-1}(A+\imath\omega_2 E).
\end{align*}
The goal in this paper is to constructed a reduced model $\hat{G}=(\hat{A}, \hat{B}, \hat{C})$ from a given model $G = (E, A, B, C)$ that minimizes the error
\begin{align}\label{eq:flmor:error}
\Xi_\omega = \| \text{G}-\hat{\text{G}}\|_{\mathcal{H}_{2,\omega}},
\end{align}
where $\| . \|_{\mathcal{H}_{2,\omega}}$ denotes the $\mathcal{H}_2$-norm on the prescribed frequency range $\omega$. 
The transfer-function matrix of the reduced  system is 
The $\mathcal{H}_2$-norm of the error system as defined in (\ref{eq:gensys:error}) can be measured efficiently by 
\begin{align}\label{eq:flmor:eerr}
 \Xi_\omega = : \sqrt{\mathrm{Tr}\left(C_{\Xi}P_{\Xi,\omega}C_{\Xi}^T\right)} \quad \text{or}\quad \sqrt{\mathrm{Tr}\left(B_{\Xi}^T Q_{\Xi,\omega}B_{\Xi}\right)} ,
\end{align}
where $P_{\Xi,\omega}$ and $Q_{\Xi,\omega}$ are the solutions of the Lyapunov equations 
\begin{subequations}
  \label{eq:flmor:rsm}
  \begin{align}
    & A_\Xi P_{\Xi,\omega} E_\Xi^{T} + E_\Xi P_{\Xi,\omega} A_\Xi^{T}  + B_{\Xi,\omega}B_\Xi B_\Xi^{T}+B_\Xi B_\Xi^{T}B_{\Xi,\omega}^{\ast}  = 0  \quad \text{and} \\
   & A_\Xi^{T} Q_{\Xi,\omega} E_\Xi + E_\Xi^{T} Q_{\Xi,\omega} A_\Xi  + C_{\Xi,\omega}^{\ast} C_\Xi^{T}C_\Xi + C_\Xi^{T}C_{\Xi} C_{\Xi,\omega}   = 0,
  \end{align}
\end{subequations}
where 
\begin{align*}
 & B_{\Xi,\omega} = \frac{\imath}{2\pi}\ln(A_\Xi+\imath\omega_2 E_\Xi)(A_\Xi+\imath\omega_1 E_\Xi)^{-1}\,\,\text{and}\\
 & C_{\Xi,\omega} = \frac{\imath}{2\pi}\ln(A_\Xi+\imath\omega_1 E_\Xi)^{-1}(A_\Xi+\imath\omega_2 E_\Xi).
\end{align*}

Due to the structure of $E_\Xi$, $A_\Xi$, $B_\Xi$ and $C_\Xi$ we can Partition $P_{\Xi,\omega}$, $Q_{\Xi,\omega}$, $B_{\Xi,\omega}$  and $C_{\Xi,\omega}$ as follows  
\begin{align*}
& P_{\Xi,\omega} = \begin{bmatrix} P_\omega & M_\omega\\ M_\omega^T & \hat{P}_\omega \end{bmatrix}, \quad
Q_{\Xi,\omega} = \begin{bmatrix} Q_\omega & N_\omega \\ N_\omega^T & \hat{Q}_\omega\end{bmatrix}, \\
& B_{\Xi,\omega} = \begin{bmatrix} B_\omega & 0\\ 0 & \hat{B}_\omega\end{bmatrix}, \qquad
C_{\Xi,\omega} = \begin{bmatrix} C_\omega & 0 \\ 0 & \hat{C}_\omega \end{bmatrix}. 
\end{align*}
Therefore (\ref{eq:flmor:rsm}) yields

\begin{subequations}  \label{eq:flmor:splitedlyap}
  \begin{align}
& A P_\omega E^{T} + E P_\omega A^{T} +B_\omega B B^{T} + BB^TB_\omega^\ast=0,  \label{eq:flmor:splitedlyap1} \\
& \hat{A} \hat{P}_\omega  +  \hat{P}_\omega \hat{A}^{T}  +\hat{B}_\omega \hat{B} \hat{B}^{T} + \hat{B}\hat{B}^T\hat{B}_\omega^\ast =0,  \label{eq:flmor:splitedlyap2} \\
& A M_\omega  + E M_\omega \hat{A}^{T} + B_\omega B \hat{B}^{T} + B\hat{B}^T\hat{B}_\omega^\ast=0,  \label{eq:flmor:splitedlyap3} \\
& A Q_\omega E^{T} + E Q_\omega A^{T} +C_\omega^\ast C^T C + C^T C C_\omega=0,  \label{eq:flmor:splitedlyap4} \\
& \hat{A} \hat{Q}_\omega  +  \hat{Q}_\omega \hat{A}^{T}  +\hat{C}_\omega^\ast \hat{C}^T\hat{C} + \hat{C}^T\hat{C} \hat{C}_\omega =0,  \label{eq:flmor:splitedlyap5} \\
& A N_\omega  + E N_\omega \hat{A}^{T} + C_\omega^\ast C^T\hat{C} - C^T\hat{C} \hat{C}_\omega=0,  \label{eq:flmor:splitedlyap6} 
\end{align}
\end{subequations}
with 
\begin{align*}
 & \hat{B}_\omega = \frac{\imath}{2\pi}\ln(\hat{A}+\imath\omega_2 I)(\hat{A}+\imath\omega_1 I)^{-1}\,\,\text{and}\\
 & \hat{C}_\omega = \frac{\imath}{2\pi}\ln(\hat{A}+\imath\omega_1 I)^{-1}(\hat{A}+\imath\omega_2 I).
\end{align*}
Following the discussion in the above section here we also construct reduced order model by constructing the reduced matrices as in (\ref{eq:gensys:reducedmatrices}). We solve the the Sylvesters equations (\ref{eq:flmor:splitedlyap3}) and (\ref{eq:flmor:splitedlyap6})
to construct $V=M_\omega$ and $W=N_\omega$. The constructed reduced order model is $\mathcal{H}_2$ optimal on the limited frequency range and satisfies Wilson's first-order optimality conditions. The whole procedure is summarized in Algorithm~\ref{alg:flmor:alorithm}. The main computation tasks in this algorithm is to solve sparse-dense Sylvesters equations   
(\ref{eq:flmor:alg:se1}) and (\ref{eq:flmor:alg:se2}). Following section will presents how to solve them efficiently.
\begin{algorithm}[t]
	\SetAlgoLined
	\SetKwInOut{Input}{Input}
	\SetKwInOut{Output}{Output}
	\caption{Two-sided iteration algorithm (TSIA).}
	\label{alg:flmor:alorithm}
	\Input {$E, A, B, C$.}
	\Output{$\hat{A}, \hat{B}, \hat{C}$, $\hat{D}:= D$.}
	Choose matrices $W_0\in \mathbb{R}^{n\times r}$ and $V_0\in\mathbb{R}^{n\times r}$  such that $W_0^TV_0=I$.\\
	Construct the reduced-order matrices
	 $\hat{A} = W_0^T E^{-1} A V_0, \hat{B} = W_0^TE^{-1}B$ and $\hat{C} = C V_0$.\\
	\While{($i\leq N-1$)}{%
	Compute $V_i=M_\omega$ and $W_i=N_\omega$ by solving the Sylvester 
	\begin{subequations}
  \label{eq:flmor:alg:se}
  \begin{align}
& A M_\omega  + E M_\omega \hat{A}^{T} + B_\omega B \hat{B}^{T} + B\hat{B}^T\hat{B}_\omega^\ast=0,  \label{eq:flmor:alg:se1}   \\
& A N_\omega  + E N_\omega \hat{A}^{T} + C_\omega^\ast C^T\hat{C} - C^T\hat{C} \hat{C}_\omega=0,  \label{eq:flmor:alg:se2}
\end{align}
\end{subequations}
	Compute $W_{i+1}=W_i(V_i^TW_i)^{-1}$ and $V_{i+1}=V_i$.\\
	Construct the reduced-order matrices
	 $\hat{A} = W_{i+1}^T E^{-1} A V_{i+1}, \hat{B} = W_{i+1}^TE^{-1}B$ and $\hat{C} = C V_{i+1}$.\\
	$i=i+1$. \\}
\end{algorithm}
%
\section{Solution of semi-generalized Sylvester equations}\label{sec:sgse}
%
Above section shows that to perform the frequency limited model reduction of system (\ref{eq:standard:ltisystem}) we need to solve two frequency limited matrix equations namely Sylvester equations (\ref{eq:flmor:alg:se1}) and (\ref{eq:flmor:alg:se2}).  This section discusses how to solve them efficiently. Since the Sylvester equations (\ref{eq:flmor:alg:se1}) and (\ref{eq:flmor:alg:se2}) are duel of each other we only interested to elaborate the solution of (\ref{eq:flmor:alg:se1}). Another one can be solved by applying the same procedure. For our convenient we rewrite the equation (\ref{eq:flmor:alg:se1}) as
\begin{align}\label{eq:sgse:gsylv_sol}
     AX+EX\hat{A}^T + F=0,  
\end{align}
where $F = B_\omega BB^T+BB^TB_\omega^\ast$ and $X=M_\omega$. The technique that we have followed here was presented in \cite{morSorA02} where $E=I$ is an identity matrix.  In \cite{benner2011sparse} authors generalized the idea of \cite{morSorA02} for the equation like (\ref{eq:sgse:gsylv_sol}) where $F=B\hat{B}$.\\

Considering   the \emph{Schur decomposition} of $\hat{A}$ as
$Q S Q^\ast$ such that $QQ^\ast = Q^\ast Q=I$ and 
inserting this into (\ref{eq:sgse:gsylv_sol}) we get
\begin{align}\label{eq:gsylv_sol1}
     AX + EXQ S Q^\ast + F=0.  
  \end{align}
By multiplying this from the right with $Q$ we obtain 
\begin{align}\label{eq:refgsylv_sol}
     A\underbrace{XQ}_{\tilde{X}} + E\underbrace{XQ}_{\tilde{X}}S + \underbrace{FQ}_{\tilde{F}}=0,  
  \end{align}
 Observing that  $S$ is a upper triangular matrices leads to a formula for the 
 first column of $\tilde{X}$:
 \begin{align}\label{eq:refgsylv_sol1}
    & A\tilde{X}_1 + E\tilde{X}_1S_{1,1} + \tilde{F}_1 =0, \\
     \Leftrightarrow\quad & (A+S_{11}E)\tilde{X}_1 =- \tilde{F}_1.
  \end{align}
For all other columns we have to take care of the linear combination of  $E$ matrix. If we consider the second column of the solution
 \begin{align}\label{eq:refgsylv_sol2}
   (A+S_{22}E)\tilde{X}_2 = - \tilde{F}_2- S_{12}E\tilde{X}_1.
  \end{align}
In this way the arbitrary column $j$ of $\tilde{X}$ we find
  \begin{align}\label{eq:refgsylv_sol3}
  (A+S_{jj}E)\tilde{X}_j = -\tilde{F}_j-E\sum_{i=1}^{j-1}S_{ij}\tilde{X}_i.
  \end{align}
  To obtain the solution of the original system we multiply $\tilde{X}$ by $U^\ast$ from the right. 
  \begin{algorithm}[t]
	\SetAlgoLined
	\SetKwInOut{Input}{Input}
	\SetKwInOut{Output}{Output}
	\caption{Solution of semi-generalized Sylvester equations.}
	\label{alg:sgse:SE}
	\Input {$E, A, \hat{A}, F$ from (\ref{eq:sgse:gsylv_sol}).}
	\Output{$X\in \mathbb{R}^{n\times r}$ solution of (\ref{eq:sgse:gsylv_sol}).}
	Compute the Schur decomposition $\hat{A}=Q S Q^\ast$ and  Define $\tilde{F}=FQ$\\
		\For{$i= 1, \cdots, r$}{%
	Compute $\hat{F} = -\tilde{F}_j-E\sum_{i=1}^{j-1}S_{i,j}\tilde{X}_i$\\
	Solve 
		$(A+S_{jj}E)\tilde{X}_j= \hat{F}$
}
 $X=\tilde{X}Z^\ast$.
\end{algorithm}
\section{FLMOR of structured index-1 systems}
\label{sec:flmor-1:index1sys}
The index 1 descriptor system that we consider in the section has the following form. 
\begin{equation}
\label{eq:flmor-1:index1sys}
\begin{array}{rcl}
  E_1 \dot{x}(t) &= J_1 x(t)+J_2 z(t)+B_1 u(t)\\
           0 &= J_3 x(t)+J_4 z(t)+B_2 u(t)\\
         y(t) &= C_1 x(t)+C_2 z(t)+D_a u(t),
  \end{array}
\end{equation}
Where $x(t)\in \mathbb{R}^{n_1}$ is the vector with differential variables and $z(t)\in \mathbb{R}^{n_2}$ is the vector with algebraic variables. Model reduction of such descriptor system has been discussed in a couple of previous research articles, e.g., \cite{morSty02,morGugSW13,morUddSKetal12, morFreRM08,morUdd11,morUddSKetal12,morHOSU19,morUddSF17} on an unrestricted frequency limit.  In \cite{morSty02} author uses spectral projector to split the system into finite and infinite parts and balancing based method was applied to the finite part. Other papers implemented MOR without computing the spectral projector, rather  eliminating the algebraic part the system was converted into an ODE system. However, practical implementation was carried out without computing the ODE system explicitly.   This paper generalizes this idea for of the FLMOR  discussed in the above section. 

By eliminating the algebraic variables i.e., $z(t)\in \mathbb{R}^{n_2}$ of the system we obtain  a generalized system (\ref{eq:standard:ltisystem}) where the coefficients matrices are defined as
\begin{equation}\label{eq:flmor-1:index0}
\begin{aligned}
  & E:=E_{1},\quad A:=J_{1}-J_{2} {J_{4}}^{-1} J_{3}, \quad B:=B_{1}-J_{2} {J_{4}}^{-1} B_{2}, \\ 
  & C:=C_{1}-C_{2} {J_{4}}^{-1} J_{3}, \quad D:=D_a-C_{2} {J_{4}}^{-1} B_{2}.
 \end{aligned}
 \end{equation}
The index-1 and generalized systems are equivalent since the responses of the systems are same and their finite eigenvalues are coincided.  Such structured system are arising in power system model \cite{morFre08}.
For the FLMOR of index-1 system we define $V$ and $W$ by solving the corresponding Sylvesters equations of the generalized system as discussed in Section~\ref{sec:flmor}.
 Now, applying these transformations the reduced system matrices can be constructed as:
\begin{equation}
\label{eq:rmat4dae}
\begin{aligned}
& \hat{A}:= \hat{J}_1-\hat{J}_2 {J}_4^{-1}\hat{J}_3,\quad
\hat{B}:= \hat{B}_1-\hat{J}_2 \hat{J}_4^{-1} B_2,\\
&\hat{C}:= \hat{C}_1- C_2 J_4^{-1}\hat{J}_3,\quad \hat{D}:= D_a-C_2{J}_4^{-1}B_2,
\end{aligned}
\end{equation}
where
$\hat{J}_1=W^TE_1^{-1} J_1 V$,
 $\hat{J}_2= W^T J_2$, $\hat{J}_3=  J_3 V$,
 $\hat{B}_1= W^TE_1^{-1} B_1$ and $\hat{C}_1= C_1 V$.
To compute the the projectors $V$ and $W$ by solving the corresponding   Sylvesters equations is a challenging task since the input matrices in (\ref{eq:rmat4dae}) are highly dense. In the following text we discuss how to solve the Sylvesters equations related to the index-1 system (\ref{eq:flmor-1:index1sys}) efficiently.

To solve the Sylvesters equations we can use Algorithm~\ref{alg:sgse:SE}. In this algorithm the main expensive task is to solve a linear system at each iteration step. At Step~4 of the algorithm we need to solve the linear system 
\begin{align*}
    (A+S_{ii}E)\tilde{X}_j=\hat{F}.
\end{align*}
Plugging $A$ and $E$ from (\ref{eq:flmor-1:index0}) we obtain 
 \begin{align*}
    (J_1-J_2J_4^{-1}J_3+S_{ii}E_1)\tilde{X}_j=\hat{F},
\end{align*}
which can be rewritten as 
 \begin{align}\label{eq:flmor-1:dansels}
    (J_1+S_{ii}E_1-J_2J_4^{-1}J_3)\tilde{X}_j=\hat{F}.
\end{align}
A close observation reveal that instead of this we can solve the following linear system 
\begin{align}\label{eq:flmor-1:sparsels}
\begin{bmatrix}
    J_1+S_{ii}E_1 & J_2\\J_3 & J_4
\end{bmatrix}
\begin{bmatrix}
 \tilde{X}_j \\ \Gamma 
\end{bmatrix}=
\begin{bmatrix}
 \hat{F} \\ 0 
\end{bmatrix}
 \end{align}
 for $\tilde{X}_j$.  Although, linear system in (\ref{eq:flmor-1:sparsels}) is larger than the system in (\ref{eq:flmor-1:dansels}) it is sparse and hence can be solved by any sparse solvers (e.g., direct \cite{Dav06,DufER89} or iterative e.g., \cite{Vor03,Saa03a})  efficiently.
  \begin{algorithm}[t]
	\SetAlgoLined
	\SetKwInOut{Input}{Input}
	\SetKwInOut{Output}{Output}
	\caption{Sylvester equation for index-1 system.}
	\label{alg:sgse:SE-index1}
	\Input {$E_1, J_1, J_2, J_3, J_4, B_1, B_2, \hat{A}$.}
	\Output{$X\in \mathbb{R}^{n\times r}$.}
	Form $F=B_\omega B\hat{B}^T+B\hat{B}^T\hat{B}_\omega^\ast$, where
	\begin{align*}
 & B_\omega = \frac{\imath}{2\pi}\ln(J_1+\imath\omega_2 E_1-J_2J_4^{-1}J_3)(J_1+\imath\omega_1 E_1-J_2J_4^{-1}J_3)^{-1},\\
 & \hat{B}_\omega = \frac{\imath}{2\pi}\ln(\hat{A}+\imath\omega_1 I)^{-1}(\hat{A}+\imath\omega_2 I),\quad \text{and}\quad B= B_1-J_2J_4^{-1}B_2
\end{align*}
	Compute the Schur decomposition $\hat{A}=Q S Q^\ast$ and  Define $\tilde{F}=FQ$\\
		\For{$i= 1, \cdots, r$}{%
	Compute $\hat{F} = -\tilde{F}_j-E_1\sum_{i=1}^{j-1}S_{i,j}\tilde{X}_i$\\
Solve 
\begin{align*}
\begin{bmatrix}
    J_1+S_{ii}E_1 & J_2\\J_3 & J_4
\end{bmatrix}
\begin{bmatrix}
 \tilde{X}_j \\ \Gamma 
\end{bmatrix}=
\begin{bmatrix}
 \hat{F} \\ 0 
\end{bmatrix}
 \end{align*}
 for $\tilde{X}_j$
}
 $X=\tilde{X}Z^\ast$.
\end{algorithm}
\section{Numerical results}\label{sec:numreslt}
\newlength\figwidth
\setlength{\figwidth}{.35\linewidth}
\newlength\figheight
\setlength{\figheight}{.5\linewidth}
\pgfplotsset{every axis plot/.append style={line width=1.5pt}}
\tikzset{mark options={solid,mark size=3,line width=.5pt,mark repeat=20}}
To asses the efficiency of our proposed techniques in this section we discuss the numerical results. For our convenient we have splitted the section into several subsections. 
\subsection{Model examples}
 We consider the following model examples for the numerical experiments.
 
\begin{example}[International Space Station (ISS)]\label{exm:num:issm}
This is a model of stage 1R (Russian Service Module) of the ISS. It has n=270 states, p=3 inputs and m=3 outputs. The details of the model can be obtain in \cite{morGugAB01}. 
\end{example}
\begin{example}[Clamped beam model (CBM)]
This structural model was  obtained by spatial discretization of an appropriate partial differential equation (see \cite{morAntSG01}). The dimension of the model is n=348 and it is a single input single output i.e., SISO system. The input represents the force applied to the structure at the free end, and the output is the resulting displacement.
\end{example}
\begin{example}[Triple chain oscillator (TCO) model]
Although this example was originated in \cite{TruV09}, the setup was described in \cite{Saa09} which resulted in second-order system. We convert into firs-order form in which the dimension of the system is n=10000. It is also an SISO system and input, output matrices are transpose of each other.
\end{example}
\begin{example}[Power system model]\label{me:psm}\index{power system}
We consider  several Brazilian Power System (BPS) models from \cite{morFreRM08} which are in index-1 descriptor form. 
Table~\ref{tab:me:bipsm} shows number of differential ($n_1$) and algebraic $n_2$ variables and inputs/outputs of several models which are used for the numerical experiments here. 
\begin{table}
\centering
       \begin{tabular}{|c|c|c|c|}
        \hline
     \textbf{model name} & \textbf{$n_1$}  & \textbf{$n_2$}& \textbf{$m/p$} \\
        \hline
                  BPS-606 & 606   &1142  &\\
               
            BPS-1142 & 1142 & 8593 & 4/4\\
                 
                  BPS-1450 & 1450 & 9855 & \\
            BPS-1693 & 1693 & b11582 & \\
        \hline
    \end{tabular}
    \caption{Dimension of differential and algebraic variables of different Brazilian power system models.}
    \label{tab:me:bipsm}
\end{table}
\end{example}
\subsection{Setup Hardware and Software}
The experiments {are} carried out {with}  Python 3.7.9 on a board with AMD $\text{Ryzen}^{\text{TM}}$ $\text{Threadripper}^{\text{TM}}$ 1920X 12-Core Processor with a 3.5~GHz clock speed and 128~GB RAM. 

\subsection{Error analysis of reduce-order model}
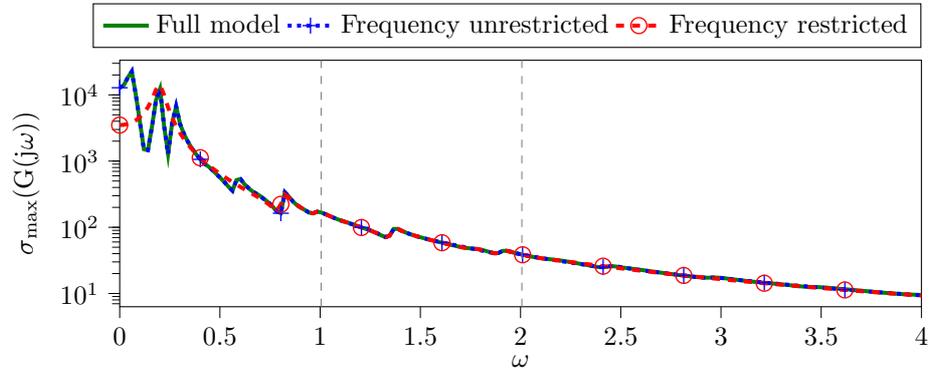
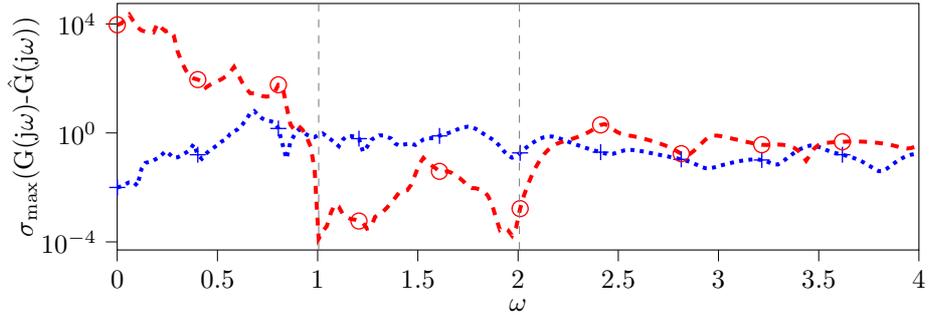
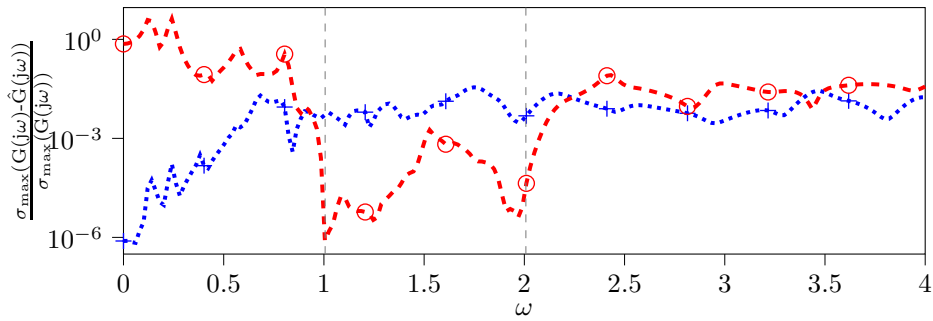
\begin{figure}[]
 \setlength{\figheight}{.4\linewidth}
   \begin{subfigure}{\linewidth}
     \centering
     \setlength{\figwidth}{\linewidth}
     \input{figures/tchain_tf}
     \caption{Sigma plot}
     \label{fig:sigmaplot}
   \end{subfigure}
   \begin{subfigure}{\linewidth}
     \centering
     \setlength{\figwidth}{\linewidth}
     \input{figures/tchain_abs_err}
     \caption{Absolute error}
     \label{fig:abserror}
   \end{subfigure}
   \begin{subfigure}{\linewidth}
     \centering
     \setlength{\figwidth}{\linewidth}
     \input{figures/tchain_rel_err}
     \caption{Relative error}
     \label{fig:relatverr}
   \end{subfigure}
   \caption{Comparison of original and 30 dimensional reduced systems
     on the frequency range [1,2] for TCO.}
   \label{fig:tco}
 \end{figure}

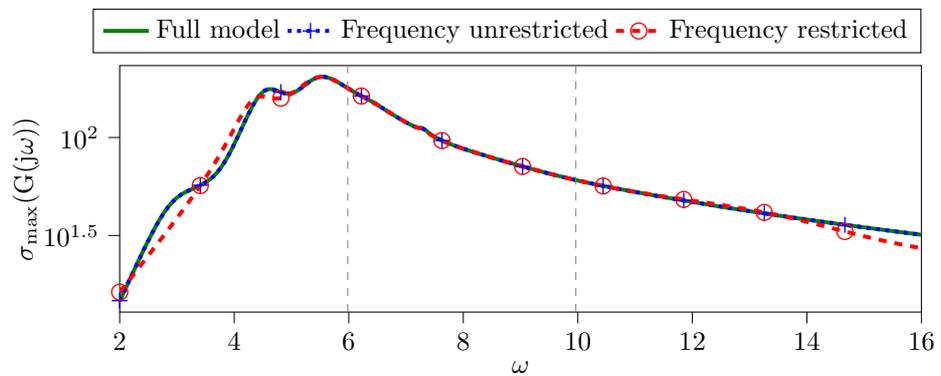
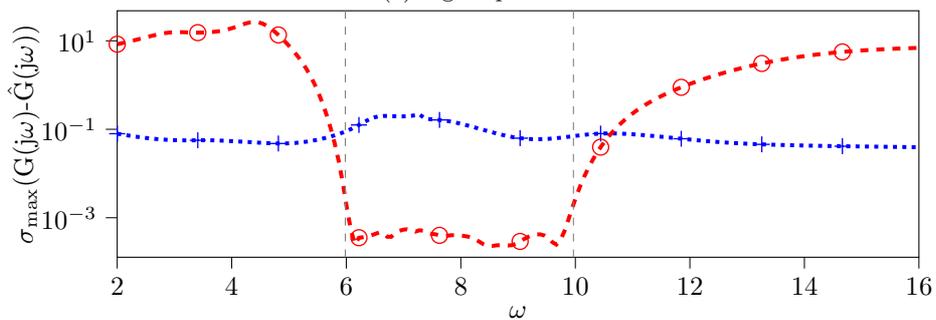
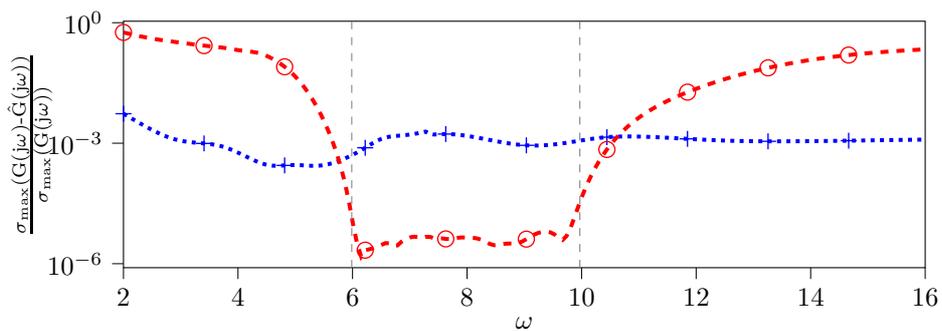
\begin{figure}[]
 \setlength{\figheight}{.4\linewidth}
   \begin{subfigure}{\linewidth}
     \centering
     \setlength{\figwidth}{\linewidth}
     \input{figures/fl_tf_1693_new}
     \caption{Sigma plot}
     \label{fig:sigmaplot}
   \end{subfigure}
   \begin{subfigure}{\linewidth}
     \centering
     \setlength{\figwidth}{\linewidth}
     \input{figures/fl_abs_err_1693_new}
     \caption{Absolute error}
     \label{fig:abserror}
   \end{subfigure}
   \begin{subfigure}{\linewidth}
     \centering
     \setlength{\figwidth}{\linewidth}
     \input{figures/fl_rel_err_1693_new}
     \caption{Relative error}
     \label{fig:relatverr}
   \end{subfigure}
   \caption{Comparison of original and 45 dimensional reduced systems
     on the frequency range [6,10] for BPS-1693.}
   \label{fig:bps}
 \end{figure}

The proposed techniques was applied to the all model examples mentioned above.  For the ISS, CBM, TCO  we apply Algorithm~\ref{alg:flmor:alorithm} to obtain frequency limited reduced order model. On the other for the BPS models we have applied the techniques discussed in Section \ref{sec:flmor-1:index1sys}. We have computed different dimensional reduced order models for different model examples which are mentioned in Table \ref{tab:roms_table}. The table also shows $\mathcal{H}_2$ norm of the error systems in both frequency restricted and unrestricted cases. For all the models frequency restricted reduced order model show much more better accuracy than the frequency frequency underused ones within the assigned frequency intervals. 
\begin{table}[]
\begin{tabularx}{1.0\textwidth} { 
  | >{\raggedright\arraybackslash}X
  | >{\centering\arraybackslash}X 
  | >{\centering\arraybackslash}X 
  | >{\centering\arraybackslash}X | }
 \hline
 Model & $r$ & $\Xi_\omega$ & $\Xi$ \\
 \hline
 ISS  & 30 & 7.3244$\times 10^{-8}$ & 2.8000$\times 10^{-3}$  \\
 \hline
 TCO  & 30 & 3.4850$\times 10^{-4}$ & 9.6000$\times 10^{-3}$  \\
 \hline
 BPS-606  & 30 & 7.8992$\times 10^{-4}$ & 1.0600$\times 10^{-2}$  \\
 \hline
 BPS-1142 & 35 & 3.4421$\times 10^{-5}$ & 5.3000$\times 10^{-3}$ \\
 \hline
 BPS-1450 & 35 & 9.3818$\times 10^{-8}$ & 1.5063$\times 10^{-5}$  \\
 \hline
 BPS-1693 & 45 & 4.2995$\times 10^{-6}$ & 1.3000 $\times 10^{-3}$  \\
 \hline
\end{tabularx}
\caption{Dimension of reduced order models and $\mathcal{H}_2$ norms of the error systems with and without limited frequency intervals.}
\label{tab:roms_table}
\end{table}
We also have investigated the frequency domain analysis including the errors of the original and reduced order models by using sigma plots. Exemplary, we have depicted the frequency responses of TCO and BPS models only. Figures \ref{fig:tco} and \ref{fig:bps} show  comparisons  of the frequency responses of the of the original and the reduced order models of TCO and BPS models, respectively . In both the figures from absolute and the relative errors we observe that frequency restricted reduced order models approximate with the original models with  higher accuracy within on the prescribed frequency ranges.  
%
%
%

%

\subsection{Comparison of sparse and dense system}
\begin{figure}
    \centering
    \input{figures/compare}
    \caption{Comparison of computational time to solve the Sylvester equation for different dimensional BPS model using Algorithms \ref{alg:sgse:SE} and \ref{alg:sgse:SE-index1}.}
    \label{fig:time_comp}
\end{figure}
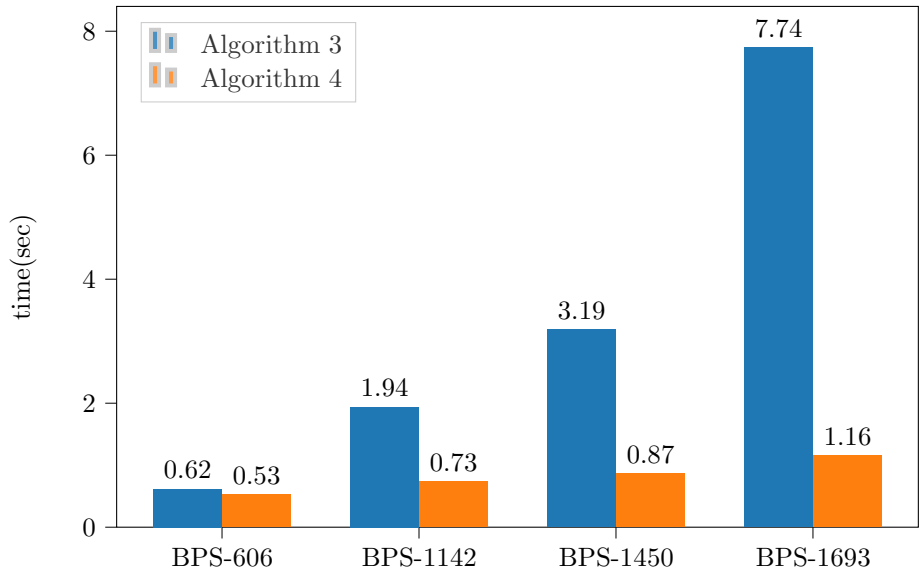

We already mentioned  that in the TSIA the main computational task is solving two Sylvester equations. This paper presents Algorithm~\ref{alg:sgse:SE-index1} to solve the Sylvester equations efficiently for index-1 system. We know that by converting index-1 system into a generalized system we can solve the Sylvester equation using Algorithm~\ref{alg:sgse:SE}. Figure \ref{fig:time_comp} shows the time comparisons of Algorithms~\ref{alg:sgse:SE} and \ref{alg:sgse:SE-index1} for solving Sylvester equations of index-1 system. We see that if the dimension of the system is increased then the computational times for Algorithms~\ref{alg:sgse:SE} is increased rapidly. On the other hand the computational time with Algorithm~\ref{alg:sgse:SE-index1} is nominal in compare to Algorithms~\ref{alg:sgse:SE}.

\section{Conclusions}
In this paper we have discussed the frequency limited $\mathcal{H}_2$ optimal model order reduction of large-scale sparse dynamical systems. For this purpose we have solved two mixed (Sparse and dense) Sylvester equations which are dual of each other. We have shown that how to solve the Sylvester equations efficiently without loosing the sparsity of large sparse system matrices. The ideas are also generalized to index-1 descriptor systems. Index-1 system can be converted into a generalized system by eliminating the algebraic equations, which however convert the system from sparse to dense. We have discussed how to perform the model reduction without converting the dense system explicitly.  Numerical experiments  has been carried out  to demonstrate the approximation accuracy and computational efficiency of the proposed algorithm using Python Programming Language.

\begin{paragraph}{Acknowledgments.}
This research work was funded by NSU-CTRG research grant under the project No.: CTRG-19/SEPS/05. It was also supported by National Natural Science Foundation of China under Grant No. (61873336, 61873335), the Fundamental Research Funds for the Central Universities under Grant (FRF-BD-19-002A), and the High-end foreign expert program of Shanghai University,
\end{paragraph}
%
  
\label{sec:numeric}
\bibliographystyle{IEEEtran}     
\bibliography{morbook}
\end{document}

%% file: figures/tchain_tf.tex
\begin{tikzpicture}

\begin{axis}[
width=\figwidth,
height=\figheight,
legend cell align={left},
log basis y={10},
tick align=outside,
tick pos=left,
x grid style={white!69.0196078431373!black},
xlabel={\(\displaystyle \omega\)},
xmin=0, xmax=4,
xtick style={color=black},
y grid style={white!69.0196078431373!black},
ylabel={$\sigma{}_{\text{max}}\text{(G(j}\omega\text{))}$},
ymin=6.30890748742499, ymax=33362.8570156218,
ymode=log,
ytick style={color=black},
legend columns=3,
legend style={nodes=right,anchor=south east,at={(1.00,1.05)}}
]
\addplot [green!50.1960784313725!black]
table {%
0 12821.6325067803
0.0201005025125628 13913.7154698806
0.0402010050251256 18275.4943209409
0.0603015075376884 22595.4858110759
0.0804020100502513 9986.32364113046
0.100502512562814 3954.13649209261
0.120603015075377 1511.61023289094
0.14070351758794 1462.81237177828
0.160804020100503 3304.37429026273
0.180904522613065 8315.60574171912
0.201005025125628 12285.2322341236
0.221105527638191 3451.01629850161
0.241206030150754 1279.44524520043
0.261306532663317 3743.43875320035
0.281407035175879 6707.87575122935
0.301507537688442 3596.11888089392
0.321608040201005 2522.54036246287
0.341708542713568 1917.98898620112
0.361809045226131 1530.93126630967
0.381909547738693 1262.72740522316
0.402010050251256 1063.01162745949
0.422110552763819 913.755543966173
0.442211055276382 831.501027552212
0.462311557788945 732.685702924804
0.482412060301508 638.977970052786
0.50251256281407 556.989655934476
0.522613065326633 481.020515577469
0.542713567839196 410.390596995868
0.562814070351759 357.969735009675
0.582914572864322 516.423389471657
0.603015075376884 535.059973201336
0.623115577889447 446.059781729224
0.64321608040201 388.047692126102
0.663316582914573 346.766156106278
0.683417085427136 321.429098247595
0.703517587939699 292.715550250335
0.723618090452261 263.331741698763
0.743718592964824 234.701235019443
0.763819095477387 204.716067825839
0.78391959798995 173.075665903212
0.804020100502513 164.014930085143
0.824120603015075 336.816148563751
0.844221105527638 290.580984794893
0.864321608040201 244.537124173443
0.884422110552764 216.190195487636
0.904522613065327 197.128943540772
0.924623115577889 181.796666257747
0.944723618090452 167.040921761175
0.964824120603015 162.594036907016
0.984924623115578 172.004756448106
1.00502512562814 167.522630347247
1.0251256281407 158.545856836206
1.04522613065327 149.138283532695
1.06532663316583 140.417682324597
1.08542713567839 132.680946641339
1.10552763819095 125.872272235956
1.12562814070352 120.593637086988
1.14572864321608 115.339085538815
1.16582914572864 109.688206506444
1.18592964824121 104.292726890735
1.20603015075377 99.2690882164791
1.22613065326633 95.0671611828102
1.24623115577889 90.6600503618492
1.26633165829146 85.5766064835564
1.28643216080402 80.3051924505567
1.30653266331658 75.1289670898758
1.32663316582915 71.3774349158413
1.34673366834171 74.4726839650578
1.36683417085427 94.3512414095499
1.38693467336683 94.5468970030762
1.4070351758794 88.4719548701347
1.42713567839196 83.5062520472181
1.44723618090452 79.257995443185
1.46733668341709 75.5883574883813
1.48743718592965 72.4335029295038
1.50753768844221 69.753627618349
1.52763819095477 67.184063259776
1.54773869346734 64.7281446342346
1.5678391959799 62.4895102119876
1.58793969849246 60.488035924582
1.60804020100503 58.4836480989678
1.62814070351759 56.5187107712831
1.64824120603015 54.6978047755022
1.66834170854271 52.9899915750928
1.68844221105528 51.3014512704062
1.70854271356784 49.8127395008015
1.7286432160804 48.7142411201935
1.74874371859296 48.0154911874943
1.76884422110553 47.4921772750896
1.78894472361809 46.6649928376066
1.80904522613065 45.3583993791509
1.82914572864322 43.7566578793585
1.84924623115578 42.0106248336023
1.86934673366834 40.4253135876215
1.8894472361809 40.7542217313192
1.90954773869347 43.9487485342188
1.92964824120603 44.2148027128791
1.94974874371859 42.6969977652595
1.96984924623116 41.1114020103208
1.98994974874372 39.7005698631083
2.01005025125628 38.4444192633066
2.03015075376884 37.3073275404779
2.05025125628141 36.2681526651099
2.07035175879397 35.3227550181828
2.09045226130653 34.4858655792083
2.1105527638191 33.7856464663298
2.13065326633166 33.2325644455045
2.15075376884422 32.7730617197242
2.17085427135678 32.3103660300404
2.19095477386935 31.7586449096833
2.21105527638191 31.1110202223624
2.23115577889447 30.4162268546705
2.25125628140704 29.7075158019395
2.2713567839196 28.9997387096235
2.29145728643216 28.2972476608558
2.31155778894472 27.5989503717024
2.33165829145729 26.9016756679728
2.35175879396985 26.2056494113048
2.37185929648241 25.5318469913574
2.39195979899498 24.9817907803483
2.41206030150754 24.8753286268815
2.4321608040201 25.4871060602031
2.45226130653266 25.9109598618968
2.47236180904523 25.6193534352707
2.49246231155779 25.0922706943877
2.51256281407035 24.5921465899359
2.53266331658291 24.1541196502064
2.55276381909548 23.744030213237
2.57286432160804 23.3297050504238
2.5929648241206 22.902206215911
2.61306532663317 22.4673092301469
2.63316582914573 22.0334219356
2.65326633165829 21.6065364975378
2.67336683417085 21.1898585438113
2.69346733668342 20.7847272621076
2.71356783919598 20.3914432266052
2.73366834170854 20.0100071548988
2.75376884422111 19.6407207012636
2.77386934673367 19.2845088445375
2.79396984924623 18.943226783072
2.81407035175879 18.6198282863303
2.83417085427136 18.3179753735997
2.85427135678392 18.040587654235
2.87437185929648 17.7881131915305
2.89447236180905 17.5619364343009
2.91457286432161 17.3827095780903
2.93467336683417 17.3106488787307
2.95477386934673 17.3685057953634
2.9748743718593 17.3996718442466
2.99497487437186 17.2704297476081
3.01507537688442 17.0255629541944
3.03517587939699 16.7376470990859
3.05527638190955 16.4433299071999
3.07537688442211 16.1562634641286
3.09547738693467 15.8808610178289
3.11557788944724 15.6183895193241
3.1356783919598 15.3692839486219
3.15577889447236 15.1339651410437
3.17587939698492 14.9129840591707
3.19597989949749 14.7067672940471
3.21608040201005 14.514994506079
3.23618090452261 14.3351197341245
3.25628140703518 14.1610421273833
3.27638190954774 13.9850329648704
3.2964824120603 13.8016282721079
3.31658291457286 13.6096919432968
3.33668341708543 13.4117135790063
3.35678391959799 13.2122164210412
3.37688442211055 13.0173219557169
3.39698492462312 12.8362342080366
3.41708542713568 12.6838838673803
3.43718592964824 12.5795131473945
3.4572864321608 12.5285268478207
3.47738693467337 12.4959852828862
3.49748743718593 12.4310249075937
3.51758793969849 12.3205485778606
3.53768844221106 12.1827020979678
3.55778894472362 12.0363284782397
3.57788944723618 11.8909930046301
3.59798994974874 11.7488878904603
3.61809045226131 11.6084321416625
3.63819095477387 11.4673178094536
3.65829145728643 11.3243900865809
3.67839195979899 11.1800173621136
3.69849246231156 11.035460180379
3.71859296482412 10.8921175271352
3.73869346733668 10.7511267122748
3.75879396984925 10.6133208143592
3.77889447236181 10.4793439050469
3.79899497487437 10.3497800649554
3.81909547738693 10.2252716217061
3.8391959798995 10.1066333848172
3.85929648241206 9.99497600429026
3.87939698492462 9.89186826296484
3.89949748743719 9.79952200401358
3.91959798994975 9.7206200660905
3.93969849246231 9.65656327809039
3.95979899497487 9.60340730013148
3.97989949748744 9.5499978318609
4 9.48461632628172
};
\addlegendentry{Full model}
\addplot [blue, dotted,mark=+]
table {%
0 12821.642474158
0.0201005025125628 13913.7233597955
0.0402010050251256 18275.4872012656
0.0603015075376884 22595.4977436198
0.0804020100502513 9986.30843597125
0.100502512562814 3954.14500337793
0.120603015075377 1511.56346484502
0.14070351758794 1462.87443551706
0.160804020100503 3304.28919828018
0.180904522613065 8315.61246026633
0.201005025125628 12285.3482654339
0.221105527638191 3450.90124532565
0.241206030150754 1279.61369657393
0.261306532663317 3743.48797984567
0.281407035175879 6707.89389529496
0.301507537688442 3595.99530701499
0.321608040201005 2522.55080480967
0.341708542713568 1917.89747395901
0.361809045226131 1531.13085906617
0.381909547738693 1262.9722084373
0.402010050251256 1062.85388408448
0.422110552763819 913.751522766083
0.442211055276382 831.364855215107
0.462311557788945 732.812475613312
0.482412060301508 639.219909590923
0.50251256281407 556.841654379528
0.522613065326633 480.627253311336
0.542713567839196 410.895542572801
0.562814070351759 357.517917797296
0.582914572864322 515.792979330605
0.603015075376884 533.740951734857
0.623115577889447 446.080993763205
0.64321608040201 390.909653210518
0.663316582914573 350.401322016112
0.683417085427136 317.292921013186
0.703517587939699 288.346976731612
0.723618090452261 261.604648527226
0.743718592964824 235.323677592951
0.763819095477387 207.20987394436
0.78391959798995 174.242385694034
0.804020100502513 163.19523836041
0.824120603015075 336.417466943114
0.844221105527638 290.488315280755
0.864321608040201 244.697747863781
0.884422110552764 217.039593464664
0.904522613065327 197.005018306133
0.924623115577889 180.699484310147
0.944723618090452 167.430656011262
0.964824120603015 162.839675607073
0.984924623115578 171.722342409041
1.00502512562814 168.318740883162
1.0251256281407 158.030015604135
1.04522613065327 148.398960404148
1.06532663316583 140.06985530517
1.08542713567839 132.756727496012
1.10552763819095 126.197452542633
1.12562814070352 120.210587835744
1.14572864321608 114.668782956062
1.16582914572864 109.476068056944
1.18592964824121 104.552890944079
1.20603015075377 99.8253325349706
1.22613065326633 95.2153250515054
1.24623115577889 90.6291954236853
1.26633165829146 85.9427574240585
1.28643216080402 80.9904550528561
1.30653266331658 75.6466399229422
1.32663316582915 70.772137086557
1.34673366834171 74.6036246823607
1.36683417085427 93.897656725693
1.38693467336683 94.8284865057169
1.4070351758794 88.5805677500579
1.42713567839196 83.1805219172453
1.44723618090452 78.8879338398046
1.46733668341709 75.3425813876971
1.48743718592965 72.2951509803722
1.50753768844221 69.5964924389663
1.52763819095477 67.1555105092607
1.54773869346734 64.9138548040672
1.5678391959799 62.8322208343802
1.58793969849246 60.8828367031793
1.60804020100503 59.0451810651623
1.62814070351759 57.3034306356447
1.64824120603015 55.6448628310324
1.66834170854271 54.0587986872266
1.68844221105528 52.5358524114029
1.70854271356784 51.0673477860656
1.7286432160804 49.6448137233388
1.74874371859296 48.2595122834414
1.76884422110553 46.9020293985247
1.78894472361809 45.5622155530675
1.80904522613065 44.2307614704467
1.82914572864322 42.9078531693795
1.84924623115578 41.6428544238633
1.86934673366834 40.7087674071391
1.8894472361809 41.1377447165003
1.90954773869347 43.6317604107377
1.92964824120603 44.2589602205084
1.94974874371859 42.8362746567993
1.96984924623116 41.2338778110624
1.98994974874372 39.8378173081659
2.01005025125628 38.629581487045
2.03015075376884 37.557293972415
2.05025125628141 36.5830244263182
2.07035175879397 35.6819657326413
2.09045226130653 34.8378022002329
2.1105527638191 34.0394329606066
2.13065326633166 33.2789799324138
2.15075376884422 32.5505937931892
2.17085427135678 31.8497210297091
2.19095477386935 31.1726349734247
2.21105527638191 30.5161183180687
2.23115577889447 29.8772322574461
2.25125628140704 29.253135749632
2.2713567839196 28.6409414961545
2.29145728643216 28.0376318531751
2.31155778894472 27.440154778662
2.33165829145729 26.8461267865475
2.35175879396985 26.2566088583666
2.37185929648241 25.6860218182592
2.39195979899498 25.195542931129
2.41206030150754 24.9785431937333
2.4321608040201 25.3310295000619
2.45226130653266 25.8605225164326
2.47236180904523 25.7688953125008
2.49246231155779 25.2629914048354
2.51256281407035 24.6787399203936
2.53266331658291 24.1209204167913
2.55276381909548 23.6049508476015
2.57286432160804 23.1262608325554
2.5929648241206 22.6778807993287
2.61306532663317 22.253936934321
2.63316582914573 21.8499560981156
2.65326633165829 21.4625719121568
2.67336683417085 21.0892025405474
2.69346733668342 20.7278001124262
2.71356783919598 20.3766690309646
2.73366834170854 20.0343343894749
2.75376884422111 19.6994453642828
2.77386934673367 19.3707076913404
2.79396984924623 19.0468563920795
2.81407035175879 18.7267208643761
2.83417085427136 18.4095491131339
2.85427135678392 18.0960978368349
2.87437185929648 17.7919933184376
2.89447236180905 17.5173190176653
2.91457286432161 17.327140995839
2.93467336683417 17.3087171771186
2.95477386934673 17.4184622027273
2.9748743718593 17.4110934781041
2.99497487437186 17.2253284314028
3.01507537688442 16.9611153826415
3.03517587939699 16.6830282582595
3.05527638190955 16.4125188341243
3.07537688442211 16.1542675161505
3.09547738693467 15.9080446729626
3.11557788944724 15.672503997406
3.1356783919598 15.4462563065537
3.15577889447236 15.22812315228
3.17587939698492 15.0171574198665
3.19597989949749 14.8126067929363
3.21608040201005 14.6138708959627
3.23618090452261 14.4204650183461
3.25628140703518 14.2319921448716
3.27638190954774 14.048122085963
3.2964824120603 13.8685760533383
3.31658291457286 13.6931152491999
3.33668341708543 13.5215323730295
3.35678391959799 13.3536452451083
3.37688442211055 13.1892919717238
3.39698492462312 13.0283272410462
3.41708542713568 12.8706194552953
3.43718592964824 12.7160484871766
3.4572864321608 12.564503906753
3.47738693467337 12.4158835661937
3.49748743718593 12.2700924593194
3.51758793969849 12.1270417940689
3.53768844221106 11.9866482314002
3.55778894472362 11.8488332553924
3.57788944723618 11.7135226476215
3.59798994974874 11.5806460450638
3.61809045226131 11.4501365654121
3.63819095477387 11.3219304872006
3.65829145728643 11.1959669748006
3.67839195979899 11.0721878404015
3.69849246231156 10.950537336677
3.71859296482412 10.8309619750693
3.73869346733668 10.713410365591
3.75879396984925 10.5978330748093
3.77889447236181 10.4841824992782
3.79899497487437 10.372412752172
3.81909547738693 10.2624795612555
3.8391959798995 10.1543401766452
3.85929648241206 10.0479532870645
3.87939698492462 9.94327894350471
3.89949748743719 9.84027848937472
3.91959798994975 9.73891449635825
3.93969849246231 9.6391507053166
3.95979899497487 9.54095197166815
3.97989949748744 9.44428421475707
4 9.34911437079022
};
\addlegendentry{Frequency unrestricted}
\addplot [red, dashed,mark=o]
table {%
0 3500.26079619555
0.0201005025125628 3534.5319795917
0.0402010050251256 3641.36662317635
0.0603015075376884 3834.01567251569
0.0804020100502513 4139.18490945461
0.100502512562814 4606.87201939855
0.120603015075377 5332.75325434815
0.14070351758794 6511.02663218671
0.160804020100503 8550.08658281097
0.180904522613065 12012.738990426
0.201005025125628 14078.3526529445
0.221105527638191 10026.2908753867
0.241206030150754 6484.71626206675
0.261306532663317 4521.34179780521
0.281407035175879 3371.52013241701
0.301507537688442 2636.97419600054
0.321608040201005 2134.16041375514
0.341708542713568 1771.58676153737
0.361809045226131 1499.54073630983
0.381909547738693 1288.97469014598
0.402010050251256 1121.88338682659
0.422110552763819 986.556519674892
0.442211055276382 875.069994948049
0.462311557788945 781.881019769529
0.482412060301508 703.001705101801
0.50251256281407 635.492116159906
0.522613065326633 577.137754898944
0.542713567839196 526.237640671662
0.562814070351759 481.460845305925
0.582914572864322 441.746444069872
0.603015075376884 406.231368021348
0.623115577889447 374.19600658606
0.64321608040201 345.020302341705
0.663316582914573 318.144287637931
0.683417085427136 293.026844618482
0.703517587939699 269.095509643523
0.723618090452261 245.68375762005
0.743718592964824 222.004755874075
0.763819095477387 197.697146912141
0.78391959798995 179.17783285151
0.804020100502513 222.482645595929
0.824120603015075 313.273044600046
0.844221105527638 280.873159018606
0.864321608040201 243.659436720661
0.884422110552764 217.449122114198
0.904522613065327 197.349425023661
0.924623115577889 180.530149378459
0.944723618090452 166.696117184603
0.964824120603015 162.894577321162
0.984924623115578 171.98415854817
1.00502512562814 167.522751374302
1.0251256281407 158.545793308666
1.04522613065327 149.138499339292
1.06532663316583 140.416962723524
1.08542713567839 132.681122444732
1.10552763819095 125.873966072223
1.12562814070352 120.59323057731
1.14572864321608 115.338641327369
1.16582914572864 109.688197606875
1.18592964824121 104.293108747962
1.20603015075377 99.2687193668259
1.22613065326633 95.0670661339174
1.24623115577889 90.6603022708488
1.26633165829146 85.576570478288
1.28643216080402 80.3050388997818
1.30653266331658 75.1283456743064
1.32663316582915 71.3771007873712
1.34673366834171 74.4721196368604
1.36683417085427 94.3487375288482
1.38693467336683 94.5440201197155
1.4070351758794 88.4781026568331
1.42713567839196 83.5016708982116
1.44723618090452 79.2608673636979
1.46733668341709 75.5845366962317
1.48743718592965 72.4535598617825
1.50753768844221 69.6811150019639
1.52763819095477 67.1227108398274
1.54773869346734 64.7595970992364
1.5678391959799 62.5494403570749
1.58793969849246 60.4554583235087
1.60804020100503 58.4599900154025
1.62814070351759 56.5489426755302
1.64824120603015 54.7135694680583
1.66834170854271 52.9559517514879
1.68844221105528 51.3028189828395
1.70854271356784 49.8335129234607
1.7286432160804 48.7015203636592
1.74874371859296 48.0145851749333
1.76884422110553 47.5005726186829
1.78894472361809 46.6565138207457
1.80904522613065 45.3597994776355
1.82914572864322 43.7618071278906
1.84924623115578 42.0056904475428
1.86934673366834 40.4274326735839
1.8894472361809 40.7535048011119
1.90954773869347 43.9490111287219
1.92964824120603 44.2146958059246
1.94974874371859 42.6967604625733
1.96984924623116 41.1115629051713
1.98994974874372 39.7002414178053
2.01005025125628 38.444203011203
2.03015075376884 37.3106935644671
2.05025125628141 36.2812700006496
2.07035175879397 35.3504525795354
2.09045226130653 34.5190851524555
2.1105527638191 33.785104491825
2.13065326633166 33.1357198075088
2.15075376884422 32.5472841373168
2.17085427135678 31.9934242273507
2.19095477386935 31.4539333252378
2.21105527638191 30.9180167374211
2.23115577889447 30.382410365123
2.25125628140704 29.8479935481756
2.2713567839196 29.3171786722482
2.29145728643216 28.7925013020038
2.31155778894472 28.2760559641972
2.33165829145729 27.7693703964241
2.35175879396985 27.2734597117367
2.37185929648241 26.7889322322922
2.39195979899498 26.3160940575479
2.41206030150754 25.8550353484948
2.4321608040201 25.4056961287043
2.45226130653266 24.9679144763495
2.47236180904523 24.5414610402448
2.49246231155779 24.1260634411463
2.51256281407035 23.7214233800456
2.53266331658291 23.3272285577075
2.55276381909548 22.9431609270452
2.57286432160804 22.5689023612037
2.5929648241206 22.2041385020358
2.61306532663317 21.8485613273565
2.63316582914573 21.5018708159779
2.65326633165829 21.1637759777203
2.67336683417085 20.8339954372235
2.69346733668342 20.5122577053585
2.71356783919598 20.1983012333078
2.73366834170854 19.8918743170181
2.75376884422111 19.5927349003237
2.77386934673367 19.3006503112282
2.79396984924623 19.0153969559539
2.81407035175879 18.7367599882904
2.83417085427136 18.4645329666809
2.85427135678392 18.1985175078043
2.87437185929648 17.9385229427644
2.89447236180905 17.6843659800676
2.91457286432161 17.4358703781868
2.93467336683417 17.1928666295024
2.95477386934673 16.955191656688
2.9748743718593 16.7226885220898
2.99497487437186 16.4952061502801
3.01507537688442 16.2725990637052
3.03517587939699 16.0547271311684
3.05527638190955 15.8414553287689
3.07537688442211 15.6326535128389
3.09547738693467 15.4281962043739
3.11557788944724 15.2279623844267
3.1356783919598 15.0318352999288
3.15577889447236 14.8397022794025
3.17587939698492 14.6514545580421
3.19597989949749 14.4669871116545
3.21608040201005 14.2861984989722
3.23618090452261 14.108990711874
3.25628140703518 13.9352690330664
3.27638190954774 13.7649419008099
3.2964824120603 13.5979207802895
3.31658291457286 13.434120041258
3.33668341708543 13.273456841598
3.35678391959799 13.1158510164701
3.37688442211055 12.9612249727385
3.39698492462312 12.8095035883788
3.41708542713568 12.6606141165951
3.43718592964824 12.5144860943859
3.4572864321608 12.3710512553206
3.47738693467337 12.2302434462953
3.49748743718593 12.0919985480582
3.51758793969849 11.9562543993007
3.53768844221106 11.8229507241297
3.55778894472362 11.6920290627411
3.57788944723618 11.5634327051307
3.59798994974874 11.437106627685
3.61809045226131 11.3129974325044
3.63819095477387 11.1910532893232
3.65829145728643 11.0712238798919
3.67839195979899 10.9534603447033
3.69849246231156 10.8377152319438
3.71859296482412 10.7239424485624
3.73869346733668 10.6120972133543
3.75879396984925 10.5021360119616
3.77889447236181 10.3940165536993
3.79899497487437 10.2876977301216
3.81909547738693 10.1831395752441
3.8391959798995 10.0803032273469
3.85929648241206 9.9791508922838
3.87939698492462 9.87964580822951
3.89949748743719 9.78175221179783
3.91959798994975 9.68543530547024
3.93969849246231 9.59066122627474
3.95979899497487 9.49739701565984
3.97989949748744 9.40561059051032
4 9.3152707152547
};
\addlegendentry{Frequency restricted}
\end{axis}
\draw [gray,dashed] (2.65,0) -- (2.65,3.3);
  \draw [gray,dashed] (5.29,0) -- (5.29,3.3);
\end{tikzpicture}

%% file: figures/tchain_abs_err.tex
\begin{tikzpicture}

\begin{axis}[
width=\figwidth,
height=\figheight,
legend cell align={left},
legend style={fill opacity=0.8, draw opacity=1, text opacity=1, draw=white!80!black},
log basis y={10},
tick align=outside,
tick pos=left,
x grid style={white!69.0196078431373!black},
xlabel={\(\displaystyle \omega\)},
xmin=0, xmax=4,
xtick style={color=black},
y grid style={white!69.0196078431373!black},
ylabel={$\sigma{}_{\text{max}}\text{(G(j}\omega\text{)-}\hat{\text{G}}\text{(j}\omega\text{))}$},
ymin=5.0528006933339e-05, ymax=56554.7047967668,
ymode=log,
ytick style={color=black}
]
\addplot [blue, dotted,mark=+]
table {%
0 0.00996737769855827
0.0201005025125628 0.010922583609031
0.0402010050251256 0.01316422880172
0.0603015075376884 0.0154203141656752
0.0804020100502513 0.0157840679366266
0.100502512562814 0.0115416811783761
0.120603015075377 0.046957367996663
0.14070351758794 0.0793639062401228
0.160804020100503 0.0863832547917992
0.180904522613065 0.0979387865802761
0.201005025125628 0.118209725001257
0.221105527638191 0.155036102765314
0.241206030150754 0.213991121693286
0.261306532663317 0.165364796489667
0.281407035175879 0.122082076254289
0.301507537688442 0.127858733685038
0.321608040201005 0.151193047765781
0.341708542713568 0.186667044781706
0.361809045226131 0.236740598502116
0.381909547738693 0.393148105401265
0.402010050251256 0.158381862656136
0.422110552763819 0.112732157673061
0.442211055276382 0.172150041427156
0.462311557788945 0.229809895123268
0.482412060301508 0.291469592921999
0.50251256281407 0.370196465940424
0.522613065326633 0.493613819280994
0.542713567839196 0.574781446504736
0.562814070351759 0.902568187241878
0.582914572864322 1.12032861786227
0.603015075376884 1.40166084277338
0.623115577889447 1.97528384031134
0.64321608040201 3.06217467760898
0.663316582914573 5.24114614767123
0.683417085427136 6.35855653151352
0.703517587939699 4.51098317754148
0.723618090452261 3.3872358980813
0.743718592964824 2.89682597814289
0.763819095477387 2.71592974507965
0.78391959798995 2.44876196769337
0.804020100502513 1.45910495196171
0.824120603015075 0.400983282342969
0.844221105527638 0.113689909814863
0.864321608040201 0.467064885406406
0.884422110552764 0.851050801087943
0.904522613065327 1.31229858318274
0.924623115577889 1.15561079045648
0.944723618090452 0.870004334235617
0.964824120603015 0.737114586687445
0.984924623115578 0.70378539716408
1.00502512562814 0.863569584169454
1.0251256281407 0.958136139906339
1.04522613065327 0.7466232634422
1.06532663316583 0.580687865292733
1.08542713567839 0.441331789311216
1.10552763819095 0.325876035016691
1.12562814070352 0.570356679153121
1.14572864321608 0.766039721677932
1.16582914572864 0.768850241406992
1.18592964824121 0.733246182001962
1.20603015075377 0.618129341304923
1.22613065326633 0.204839211546516
1.24623115577889 0.306731259154069
1.26633165829146 0.547599258488365
1.28643216080402 0.688157240037571
1.30653266331658 0.811949950997765
1.32663316582915 0.811029523032955
1.34673366834171 0.706472960426912
1.36683417085427 0.626082973608283
1.38693467336683 0.410334016105487
1.4070351758794 0.338577360707497
1.42713567839196 0.370829451065788
1.44723618090452 0.394034176306197
1.46733668341709 0.395620495612101
1.48743718592965 0.361911915269793
1.50753768844221 0.440943138538448
1.52763819095477 0.579081361351987
1.54773869346734 0.64252718222844
1.5678391959799 0.637007037798772
1.58793969849246 0.646315729596841
1.60804020100503 0.778838926032907
1.62814070351759 0.902101377218258
1.64824120603015 0.979074126009294
1.66834170854271 1.07045798213162
1.68844221105528 1.25850857134987
1.70854271356784 1.46518013219135
1.7286432160804 1.61827052619337
1.74874371859296 1.69304678349932
1.76884422110553 1.62415988334032
1.78894472361809 1.4108720255215
1.80904522613065 1.14231605980753
1.82914572864322 0.911845730659374
1.84924623115578 0.739828250485573
1.86934673366834 0.59680308050188
1.8894472361809 0.466353908610941
1.90954773869347 0.321502967183702
1.92964824120603 0.212380420600166
1.94974874371859 0.149825637282606
1.96984924623116 0.123060233838232
1.98994974874372 0.137976858020461
2.01005025125628 0.185356261110011
2.03015075376884 0.251350865775307
2.05025125628141 0.330521179362446
2.07035175879397 0.421546519163303
2.09045226130653 0.52163103568247
2.1105527638191 0.621312660899807
2.13065326633166 0.701463828198308
2.15075376884422 0.742445975270126
2.17085427135678 0.737528644216265
2.19095477386935 0.685753269843097
2.21105527638191 0.612216082111273
2.23115577889447 0.539329592480213
2.25125628140704 0.474093988468739
2.2713567839196 0.417571188210576
2.29145728643216 0.369074295058814
2.31155778894472 0.327566720764209
2.33165829145729 0.292081649457594
2.35175879396985 0.261828848488135
2.37185929648241 0.236317008047842
2.39195979899498 0.215213042234175
2.41206030150754 0.196872969738374
2.4321608040201 0.177830310020193
2.45226130653266 0.165292267832312
2.47236180904523 0.170705553020566
2.49246231155779 0.189546872220267
2.51256281407035 0.212200675985437
2.53266331658291 0.230076234212036
2.55276381909548 0.238184049769096
2.57286432160804 0.2364418731416
2.5929648241206 0.22781330123068
2.61306532663317 0.215582799827538
2.63316582914573 0.202066665801877
2.65326633165829 0.188567601610838
2.67336683417085 0.175717968132784
2.69346733668342 0.163792882250429
2.71356783919598 0.152854199013842
2.73366834170854 0.142789960782753
2.75376884422111 0.133428620283749
2.77386934673367 0.124568608078644
2.79396984924623 0.115948530224053
2.81407035175879 0.107231707485256
2.83417085427136 0.0980270978205791
2.85427135678392 0.08800180052341
2.87437185929648 0.0771488152126076
2.89447236180905 0.0661647964499431
2.91457286432161 0.0566584009871012
2.93467336683417 0.0508038134767475
2.95477386934673 0.0500611693163072
2.9748743718593 0.053467944838691
2.99497487437186 0.0587824363314126
3.01507537688442 0.0647578941063584
3.03517587939699 0.0709796892067097
3.05527638190955 0.0772828249505729
3.07537688442211 0.0835529329563726
3.09547738693467 0.0896563251176328
3.11557788944724 0.0953920866450622
3.1356783919598 0.100445918761521
3.15577889447236 0.104347963080651
3.17587939698492 0.106442958264824
3.19597989949749 0.105859758644971
3.21608040201005 0.101519036331646
3.23618090452261 0.0927106938336468
3.25628140703518 0.0806569245655092
3.27638190954774 0.0701646503261797
3.2964824120603 0.0696044402891392
3.31658291457286 0.0841317920500546
3.33668341708543 0.110564113568574
3.35678391959799 0.144520257591416
3.37688442211055 0.184042168169967
3.39698492462312 0.228257084541234
3.41708542713568 0.274996633758428
3.43718592964824 0.317535914423879
3.4572864321608 0.34292209180031
3.47738693467337 0.339723176155238
3.49748743718593 0.311422384157322
3.51758793969849 0.273077022551333
3.53768844221106 0.237080777818231
3.55778894472362 0.208525719969222
3.57788944723618 0.187546536526436
3.59798994974874 0.171891820886093
3.61809045226131 0.15871292490113
3.63819095477387 0.145804887465295
3.65829145728643 0.1321046556957
3.67839195979899 0.11749163503706
3.69849246231156 0.102323399170887
3.71859296482412 0.0870825058953544
3.73869346733668 0.072257842104467
3.75879396984925 0.0584653552264939
3.77889447236181 0.0468347584357751
3.79899497487437 0.0396372191108146
3.81909547738693 0.039919090607947
3.8391959798995 0.0482863639442314
3.85929648241206 0.0623832047851859
3.87939698492462 0.0801051362057933
3.89949748743719 0.100285280695601
3.91959798994975 0.121796599794933
3.93969849246231 0.142501174368365
3.95979899497487 0.158854213865553
3.97989949748744 0.167380243353687
4 0.167132015019094
};
\addplot [red, dashed,mark=o]
table {%
0 9321.37171058478
0.0201005025125628 10421.1284835229
0.0402010050251256 15016.1926957676
0.0603015075376884 21935.798146673
0.0804020100502513 12662.5207361712
0.100502512562814 7687.37612848304
0.120603015075377 5839.01859568928
0.14070351758794 5220.50598254233
0.160804020100503 5272.73490513409
0.180904522613065 4884.6345230546
0.201005025125628 10386.3494351418
0.221105527638191 7174.71960999662
0.241206030150754 5373.58351060416
0.261306532663317 5949.20432470618
0.281407035175879 3869.46031318275
0.301507537688442 961.145728792075
0.321608040201005 401.709583142168
0.341708542713568 202.365844349801
0.361809045226131 126.21276247353
0.381909547738693 100.048868434223
0.402010050251256 90.9805644802002
0.422110552763819 78.9589759221088
0.442211055276382 46.2841601055182
0.462311557788945 58.137201399213
0.482412060301508 69.9554015623156
0.50251256281407 80.7057346404142
0.522613065326633 96.1378013981884
0.542713567839196 120.981654900056
0.562814070351759 167.13662848133
0.582914572864322 266.093908836605
0.603015075376884 141.173999875132
0.623115577889447 72.1963931798221
0.64321608040201 43.1738949856203
0.663316582914573 28.7576873305992
0.683417085427136 28.641507298248
0.703517587939699 25.9946302082685
0.723618090452261 22.7377235728418
0.743718592964824 21.4236907274171
0.763819095477387 22.631857710968
0.78391959798995 30.0967491419725
0.804020100502513 58.6377907789186
0.824120603015075 53.4586911573668
0.844221105527638 9.74966193201285
0.864321608040201 2.04185359652514
0.884422110552764 1.26525352621892
0.904522613065327 1.81059156582485
0.924623115577889 1.30361017719123
0.944723618090452 0.70884966832536
0.964824120603015 0.317621682452366
0.984924623115578 0.0449054333164823
1.00502512562814 0.000130270915923676
1.0251256281407 0.000242776819395946
1.04522613065327 0.000376614617963305
1.06532663316583 0.00113627673058823
1.08542713567839 0.00238971902972624
1.10552763819095 0.00195875852528507
1.12562814070352 0.00102825286631422
1.14572864321608 0.000794000545759082
1.16582914572864 0.000689798071929691
1.18592964824121 0.00065609295497448
1.20603015075377 0.000583314818544906
1.22613065326633 0.00050798492216419
1.24623115577889 0.000306460916100316
1.26633165829146 0.000346102092569666
1.28643216080402 0.000832074487938731
1.30653266331658 0.0010375043127048
1.32663316582915 0.00129772299562476
1.34673366834171 0.00183783157865958
1.36683417085427 0.00273924062910132
1.38693467336683 0.00417434378660919
1.4070351758794 0.00649309672509922
1.42713567839196 0.00738777497790547
1.44723618090452 0.00791050086314888
1.46733668341709 0.012771111549553
1.48743718592965 0.0302385706201283
1.50753768844221 0.0727909087557616
1.52763819095477 0.120694494904424
1.54773869346734 0.0835787122514277
1.5678391959799 0.0641579782037438
1.58793969849246 0.0510584318850633
1.60804020100503 0.0388790244132109
1.62814070351759 0.0363516628877552
1.64824120603015 0.0392385566016478
1.66834170854271 0.034305539292399
1.68844221105528 0.0266064871927464
1.70854271356784 0.0219652258611906
1.7286432160804 0.0163336891681133
1.74874371859296 0.0110862083197791
1.76884422110553 0.00953761198065681
1.78894472361809 0.00927752584671024
1.80904522613065 0.00837448392116876
1.82914572864322 0.00724887994595612
1.84924623115578 0.00493438609953152
1.86934673366834 0.00245402512200378
1.8894472361809 0.000765167931380258
1.90954773869347 0.000271667516764546
1.92964824120603 0.000321325208679289
1.94974874371859 0.000268820580642974
1.96984924623116 0.000169483344673778
1.98994974874372 0.000329283650057947
2.01005025125628 0.00166663535521908
2.03015075376884 0.00529289122600845
2.05025125628141 0.0134689389288412
2.07035175879397 0.0297737321017939
2.09045226130653 0.0589466958350274
2.1105527638191 0.105276564169534
2.13065326633166 0.168642385523723
2.15075376884422 0.242664444483501
2.17085427135678 0.321659203574016
2.19095477386935 0.392561084257628
2.21105527638191 0.448450220091536
2.23115577889447 0.498560393015758
2.25125628140704 0.550836767706168
2.2713567839196 0.6110331519955
2.29145728643216 0.684601227527304
2.31155778894472 0.778203374128515
2.33165829145729 0.901313959684552
2.35175879396985 1.06812300815913
2.37185929648241 1.29893114035361
2.39195979899498 1.61278285315083
2.41206030150754 1.9674834328228
2.4321608040201 2.10478550140562
2.45226130653266 1.8095432938792
2.47236180904523 1.39429209843909
2.49246231155779 1.1000354677348
2.51256281407035 0.938512172127456
2.53266331658291 0.858400338712132
2.55276381909548 0.808546917512857
2.57286432160804 0.760874844833362
2.5929648241206 0.70831371567534
2.61306532663317 0.652694296642033
2.63316582914573 0.596975477595618
2.65326633165829 0.542928291435176
2.67336683417085 0.491146706202771
2.69346733668342 0.441519878664972
2.71356783919598 0.393548961503954
2.73366834170854 0.346570164589618
2.75376884422111 0.300102817049618
2.77386934673367 0.254228752356941
2.79396984924623 0.210378418231283
2.81407035175879 0.173189828004316
2.83417085427136 0.153384647331384
2.85427135678392 0.164854093002869
2.87437185929648 0.211100298224278
2.89447236180905 0.290042129836054
2.91457286432161 0.408626474119233
2.93467336683417 0.568128195813971
2.95477386934673 0.723146830690654
2.9748743718593 0.800845902931707
2.99497487437186 0.796366645108123
3.01507537688442 0.753225324291737
3.03517587939699 0.700729186919106
3.05527638190955 0.649475356656749
3.07537688442211 0.601956381337078
3.09547738693467 0.558127264293801
3.11557788944724 0.517403975630761
3.1356783919598 0.479374048445308
3.15577889447236 0.444146071738346
3.17587939698492 0.412606402467457
3.19597989949749 0.386582071761311
3.21608040201005 0.368767766791066
3.23618090452261 0.361501332068779
3.25628140703518 0.363728770559404
3.27638190954774 0.370013467507795
3.2964824120603 0.373565033600653
3.31658291457286 0.369291721345767
3.33668341708543 0.354347599188838
3.35678391959799 0.326995993665945
3.37688442211055 0.28530877810916
3.39698492462312 0.226802809522461
3.41708542713568 0.151993725323597
3.43718592964824 0.0951546409151906
3.4572864321608 0.158685236077257
3.47738693467337 0.267976834553753
3.49748743718593 0.353816592346975
3.51758793969849 0.405240388448493
3.53768844221106 0.431814090519489
3.55778894472362 0.445501461960144
3.57788944723618 0.454555369127186
3.59798994974874 0.463064753956976
3.61809045226131 0.471798892009081
3.63819095477387 0.479578856454563
3.65829145728643 0.484823264712359
3.67839195979899 0.48654399836415
3.69849246231156 0.484523677780619
3.71859296482412 0.479021814328061
3.73869346733668 0.470425130054167
3.75879396984925 0.459041194757319
3.77889447236181 0.445039749119139
3.79899497487437 0.428459851584534
3.81909547738693 0.409239953541394
3.8391959798995 0.387274271966835
3.85929648241206 0.362522520373977
3.87939698492462 0.335270621445513
3.89949748743719 0.306851952834329
3.91959798994975 0.281428433880682
3.93969849246231 0.268293722303409
3.95979899497487 0.278407847178391
3.97989949748744 0.311135314502512
4 0.352431885375914
};
\end{axis}
\draw [gray,dashed] (2.65,0) -- (2.65,3.3);
  \draw [gray,dashed] (5.29,0) -- (5.29,3.3);
\end{tikzpicture}

%% file: figures/tchain_rel_err.tex
\begin{tikzpicture}

\begin{axis}[
width=\figwidth,
height=\figheight,
legend cell align={left},
legend style={fill opacity=0.8, draw opacity=1, text opacity=1, draw=white!80!black},
log basis y={10},
tick align=outside,
tick pos=left,
x grid style={white!69.0196078431373!black},
xlabel={\(\displaystyle \omega\)},
xmin=0, xmax=4,
xtick style={color=black},
y grid style={white!69.0196078431373!black},
ylabel={$\frac{\sigma{}_{\text{max}}\text{(G(j}\omega\text{)-}{\hat{\text{G}}}\text{(j}\omega\text{))}}{\sigma{}_{\text{max}}\text{(G(j}\omega\text{))}}$},
ymin=3.12329496212709e-07, ymax=9.17700220196981,
ymode=log,
ytick style={color=black}
]
\addplot [blue, dotted,mark=+]
table {%
0 7.77387566933253e-07
0.0201005025125628 7.85022780771642e-07
0.0402010050251256 7.20321353312774e-07
0.0603015075376884 6.82451100835213e-07
0.0804020100502513 1.58056843577721e-06
0.100502512562814 2.91888790421293e-06
0.120603015075377 3.10644681908891e-05
0.14070351758794 5.42543307475884e-05
0.160804020100503 2.61420914229819e-05
0.180904522613065 1.177770923998e-05
0.201005025125628 9.62209934240525e-06
0.221105527638191 4.49247668962412e-05
0.241206030150754 0.000167253051661279
0.261306532663317 4.41745697985022e-05
0.281407035175879 1.81998118006159e-05
0.301507537688442 3.55546459724532e-05
0.321608040201005 5.99368200468214e-05
0.341708542713568 9.73243569825863e-05
0.361809045226131 0.000154638293509272
0.381909547738693 0.000311348358937204
0.402010050251256 0.000148993537384587
0.422110552763819 0.000123372337839664
0.442211055276382 0.000207035272023576
0.462311557788945 0.000313654127828469
0.482412060301508 0.000456149674296159
0.50251256281407 0.000664637955114867
0.522613065326633 0.00102618038793711
0.542713567839196 0.0014005716766228
0.562814070351759 0.00252135334071603
0.582914572864322 0.00216939945150133
0.603015075376884 0.00261963314950855
0.623115577889447 0.00442829396690693
0.64321608040201 0.00789123280396647
0.663316582914573 0.0151143531609957
0.683417085427136 0.0197821434530348
0.703517587939699 0.0154108081162194
0.723618090452261 0.0128629988782595
0.743718592964824 0.0123426107148644
0.763819095477387 0.013266812780862
0.78391959798995 0.0141485052500839
0.804020100502513 0.00889617153270289
0.824120603015075 0.00119051085897407
0.844221105527638 0.000391250342465151
0.864321608040201 0.001909995821637
0.884422110552764 0.00393658370662149
0.904522613065327 0.0066570568462024
0.924623115577889 0.00635661156084171
0.944723618090452 0.00520833053998287
0.964824120603015 0.00453346629869941
0.984924623115578 0.00409166241502404
1.00502512562814 0.00515494284192779
1.0251256281407 0.00604327453915236
1.04522613065327 0.00500624819970199
1.06532663316583 0.00413543262984775
1.08542713567839 0.00332626349512124
1.10552763819095 0.00258894218105331
1.12562814070352 0.00472957523241216
1.14572864321608 0.00664163165590679
1.16582914572864 0.00700941574208181
1.18592964824121 0.0070306550021476
1.20603015075377 0.00622680587089658
1.22613065326633 0.00215467895536102
1.24623115577889 0.0033833122519767
1.26633165829146 0.00639893635643973
1.28643216080402 0.00856927452681549
1.30653266331658 0.0108074153345731
1.32663316582915 0.0113625478975143
1.34673366834171 0.00948633677226386
1.36683417085427 0.00663566227910716
1.38693467336683 0.00434000510976194
1.4070351758794 0.00382694562592727
1.42713567839196 0.00444073877074623
1.44723618090452 0.00497153850665647
1.46733668341709 0.005233881364242
1.48743718592965 0.00499647125477315
1.50753768844221 0.00632143665632748
1.52763819095477 0.00861932626957814
1.54773869346734 0.00992655027977749
1.5678391959799 0.0101938235015414
1.58793969849246 0.0106850176190657
1.60804020100503 0.0133172083368488
1.62814070351759 0.015961110310333
1.64824120603015 0.0178996968896236
1.66834170854271 0.0202011351636971
1.68844221105528 0.0245316368286807
1.70854271356784 0.0294137633640441
1.7286432160804 0.0332196599799345
1.74874371859296 0.0352604282832012
1.76884422110553 0.034198471759521
1.78894472361809 0.030234056403508
1.80904522613065 0.0251842233289342
1.82914572864322 0.0208390168456975
1.84924623115578 0.0176105033766082
1.86934673366834 0.014763103301804
1.8894472361809 0.011443082183865
1.90954773869347 0.00731540664766318
1.92964824120603 0.00480337822559825
1.94974874371859 0.00350904384674353
1.96984924623116 0.00299333585868315
1.98994974874372 0.00347543772032037
2.01005025125628 0.00482140879383564
2.03015075376884 0.00673730557361941
2.05025125628141 0.00911326205154112
2.07035175879397 0.0119341347793032
2.09045226130653 0.0151259371606715
2.1105527638191 0.0183898408313423
2.13065326633166 0.021107724904847
2.15075376884422 0.0226541536344556
2.17085427135678 0.0228263784919846
2.19095477386935 0.0215926489241992
2.21105527638191 0.0196784315569059
2.23115577889447 0.0177316402543006
2.25125628140704 0.0159587220833116
2.2713567839196 0.0143991362264242
2.29145728643216 0.0130427630093991
2.31155778894472 0.011868810818982
2.33165829145729 0.0108573775501028
2.35175879396985 0.0099913131088133
2.37185929648241 0.00925577409765285
2.39195979899498 0.00861479643819091
2.41206030150754 0.00791438668776514
2.4321608040201 0.00697726566524028
2.45226130653266 0.00637924139874809
2.47236180904523 0.0066631483675756
2.49246231155779 0.00755399439647616
2.51256281407035 0.00862879843406098
2.53266331658291 0.00952534133075182
2.55276381909548 0.0100313235634408
2.57286432160804 0.0101347990739945
2.5929648241206 0.00994722076480169
2.61306532663317 0.00959539914723148
2.63316582914573 0.00917091618326395
2.65326633165829 0.00872734052643403
2.67336683417085 0.00829255031455149
2.69346733668342 0.00788044414463105
2.71356783919598 0.00749599708638621
2.73366834170854 0.00713592752253445
2.75376884422111 0.00679346864675719
2.77386934673367 0.00645951675942887
2.79396984924623 0.00612084369531315
2.81407035175879 0.00575900625055603
2.83417085427136 0.00535141552607707
2.85427135678392 0.00487798968692419
2.87437185929648 0.00433709940913467
2.89447236180905 0.00376751144143274
2.91457286432161 0.00325946888386809
2.93467336683417 0.00293483010559871
2.95477386934673 0.002882295685428
2.9748743718593 0.00307292834700045
2.99497487437186 0.00340364641705304
3.01507537688442 0.00380356845060473
3.03517587939699 0.00424072086037626
3.05527638190955 0.00469994979038485
3.07537688442211 0.00517155053468171
3.09547738693467 0.00564555819844897
3.11557788944724 0.00610767752507623
3.1356783919598 0.00653549762612898
3.15577889447236 0.00689495199097928
3.17587939698492 0.00713760289976086
3.19597989949749 0.00719803043921285
3.21608040201005 0.00699408024502722
3.23618090452261 0.00646738189517529
3.25628140703518 0.00569569130858966
3.27638190954774 0.00501712441453868
3.2964824120603 0.00504320496950384
3.31658291457286 0.00618175579583874
3.33668341708543 0.00824384691167599
3.35678391959799 0.0109383810396309
3.37688442211055 0.0141382512314017
3.39698492462312 0.0177822467899756
3.41708542713568 0.0216807908865872
3.43718592964824 0.0252423055410255
3.4572864321608 0.0273713019867105
3.47738693467337 0.0271865858085239
3.49748743718593 0.0250520280083328
3.51758793969849 0.0221643558178926
3.53768844221106 0.0194604428403268
3.55778894472362 0.0173246950136175
3.57788944723618 0.0157721509426008
3.59798994974874 0.0146304758789692
3.61809045226131 0.01367221024892
3.63819095477387 0.0127148204914225
3.65829145728643 0.011665498511239
3.67839195979899 0.0105090744702429
3.69849246231156 0.00927223672582478
3.71859296482412 0.00799500241146023
3.73869346733668 0.00672095530433742
3.75879396984925 0.00550867690227487
3.77889447236181 0.00446924529437565
3.79899497487437 0.00382976438746047
3.81909547738693 0.00390396383438922
3.8391959798995 0.00477769026595644
3.85929648241206 0.0062414561834274
3.87939698492462 0.00809807956154319
3.89949748743719 0.0102336910570257
3.91959798994975 0.0125297150764908
3.93969849246231 0.0147569244113673
3.95979899497487 0.0165414429380058
3.97989949748744 0.0175267310318406
4 0.0176213785850223
};
\addplot [red, dashed,mark=o]
table {%
0 0.727003500190437
0.0201005025125628 0.748982434352765
0.0402010050251256 0.821657265848146
0.0603015075376884 0.970804448732872
0.0804020100502513 1.26798621707175
0.100502512562814 1.94413524769721
0.120603015075377 3.86278054265498
0.14070351758794 3.56881448589062
0.160804020100503 1.59568331005107
0.180904522613065 0.587405737449594
0.201005025125628 0.84543370749578
0.221105527638191 2.07901643730625
0.241206030150754 4.19993237753785
0.261306532663317 1.58923511694163
0.281407035175879 0.576853307468255
0.301507537688442 0.267273068723789
0.321608040201005 0.159248029930415
0.341708542713568 0.1055093881173
0.361809045226131 0.0824418216878982
0.381909547738693 0.0792323568969673
0.402010050251256 0.0855875534472149
0.422110552763819 0.0864114876714027
0.442211055276382 0.0556633829326349
0.462311557788945 0.0793480767635228
0.482412060301508 0.109480146172389
0.50251256281407 0.144896289869176
0.522613065326633 0.199862164470832
0.542713567839196 0.294796361772574
0.562814070351759 0.466901562158076
0.582914572864322 0.515263085021849
0.603015075376884 0.263847058172692
0.623115577889447 0.161853626211134
0.64321608040201 0.111259249472846
0.663316582914573 0.0829310670150448
0.683417085427136 0.0891067655492273
0.703517587939699 0.0888050880318365
0.723618090452261 0.0863463076124432
0.743718592964824 0.0912806902172555
0.763819095477387 0.110552424884508
0.78391959798995 0.173893591481562
0.804020100502513 0.357514957622935
0.824120603015075 0.158717719994499
0.844221105527638 0.0335523053543736
0.864321608040201 0.00834987163371115
0.884422110552764 0.00585250188319145
0.904522613065327 0.00918480834576365
0.924623115577889 0.00717070452404779
0.944723618090452 0.00424356894617015
0.964824120603015 0.00195346452117434
0.984924623115578 0.00026107088108363
1.00502512562814 7.77631748341382e-07
1.0251256281407 1.53127192498482e-06
1.04522613065327 2.52527123849284e-06
1.06532663316583 8.0921199650736e-06
1.08542713567839 1.80110188404525e-05
1.10552763819095 1.55614774444784e-05
1.12562814070352 8.5265930371808e-06
1.14572864321608 6.88405445604022e-06
1.16582914572864 6.28871684477003e-06
1.18592964824121 6.29087928309566e-06
1.20603015075377 5.8760972728273e-06
1.22613065326633 5.34343211519019e-06
1.24623115577889 3.38033030951501e-06
1.26633165829146 4.04435402140151e-06
1.28643216080402 1.03614033233409e-05
1.30653266331658 1.38096443075498e-05
1.32663316582915 1.8181138018687e-05
1.34673366834171 2.46779286150326e-05
1.36683417085427 2.90323750718988e-05
1.38693467336683 4.41510395256374e-05
1.4070351758794 7.33915819383695e-05
1.42713567839196 8.84697228864743e-05
1.44723618090452 9.98069761784401e-05
1.46733668341709 0.000168956066435443
1.48743718592965 0.000417466633493593
1.50753768844221 0.00104354298466068
1.52763819095477 0.00179647507233588
1.54773869346734 0.00129122675651703
1.5678391959799 0.00102669996910035
1.58793969849246 0.000844107948036605
1.60804020100503 0.000664784528273931
1.62814070351759 0.00064317926562163
1.64824120603015 0.000717369860869111
1.66834170854271 0.000647396579480187
1.68844221105528 0.000518630302532877
1.70854271356784 0.000440955990000051
1.7286432160804 0.00033529597900977
1.74874371859296 0.000230888158083999
1.76884422110553 0.000200824904813522
1.78894472361809 0.000198811256202179
1.80904522613065 0.00018462917642147
1.82914572864322 0.000165663473795051
1.84924623115578 0.000117455670299499
1.86934673366834 6.0705159817368e-05
1.8894472361809 1.87751820271477e-05
1.90954773869347 6.18146194886582e-06
1.92964824120603 7.26736723820532e-06
1.94974874371859 6.29600662137655e-06
1.96984924623116 4.12253867263467e-06
1.98994974874372 8.29417943352831e-06
2.01005025125628 4.33518150919192e-05
2.03015075376884 0.000141872698339642
2.05025125628141 0.000371370967063299
2.07035175879397 0.000842905149569095
2.09045226130653 0.00170930016819896
2.1105527638191 0.00311601449670204
2.13065326633166 0.00507461245731628
2.15075376884422 0.00740438737639991
2.17085427135678 0.00995529432489111
2.19095477386935 0.0123607630418115
2.21105527638191 0.0144145134709916
2.23115577889447 0.0163912636303606
2.25125628140704 0.0185420003267391
2.2713567839196 0.0210702985331633
2.29145728643216 0.0241932090262732
2.31155778894472 0.0281968467513322
2.33165829145729 0.0335040081074798
2.35175879396985 0.0407592649735425
2.37185929648241 0.0508749382993443
2.39195979899498 0.0645583364031496
2.41206030150754 0.0790937664516577
2.4321608040201 0.0825823652333814
2.45226130653266 0.0698369841767309
2.47236180904523 0.0544233913616079
2.49246231155779 0.0438396142434747
2.51256281407035 0.038163084653681
2.53266331658291 0.0355384651207852
2.55276381909548 0.0340526401900424
2.57286432160804 0.0326139933269127
2.5929648241206 0.030927750322292
2.61306532663317 0.0290508440488387
2.63316582914573 0.0270940882147347
2.65326633165829 0.0251279649330677
2.67336683417085 0.0231783853199065
2.69346733668342 0.0212425149051584
2.71356783919598 0.0192997110175351
2.73366834170854 0.0173198421123389
2.75376884422111 0.015279623472794
2.77386934673367 0.0131830556021111
2.79396984924623 0.0111057329693841
2.81407035175879 0.00930136547668715
2.83417085427136 0.008373449805618
2.85427135678392 0.00913795582286202
2.87437185929648 0.0118674924063779
2.89447236180905 0.0165153843325365
2.91457286432161 0.0235076397200054
2.93467336683417 0.0328195782719629
2.95477386934673 0.0416355234705164
2.9748743718593 0.0460264946431455
2.99497487437186 0.0461115708610792
3.01507537688442 0.0442408469146199
3.03517587939699 0.0418654535354241
3.05527638190955 0.0394978000394172
3.07537688442211 0.0372583910056672
3.09547738693467 0.035144647614963
3.11557788944724 0.0331278698735612
3.1356783919598 0.0311903957300685
3.15577889447236 0.0293476341196141
3.17587939698492 0.0276675949515097
3.19597989949749 0.0262859990936139
3.21608040201005 0.0254059873489185
3.23618090452261 0.0252178801972774
3.25628140703518 0.0256851697274496
3.27638190954774 0.0264578187578997
3.2964824120603 0.0270667363470151
3.31658291457286 0.0271344658559781
3.33668341708543 0.0264207550438228
3.35678391959799 0.0247495184188162
3.37688442211055 0.0219176247679623
3.39698492462312 0.0176689522679839
3.41708542713568 0.0119832164116928
3.43718592964824 0.00756425465757387
3.4572864321608 0.0126659133994561
3.47738693467337 0.0214450344240369
3.49748743718593 0.0284623830277131
3.51758793969849 0.0328914241023886
3.53768844221106 0.0354448534526278
3.55778894472362 0.0370130694559857
3.57788944723618 0.0382268637236765
3.59798994974874 0.0394134966878838
3.61809045226131 0.0406427746875309
3.63819095477387 0.0418213626258096
3.65829145728643 0.0428123069768552
3.67839195979899 0.0435190735940116
3.69849246231156 0.0439060691498937
3.71859296482412 0.0439787592389347
3.73869346733668 0.0437558911399558
3.75879396984925 0.0432514198700431
3.77889447236181 0.0424682836207718
3.79899497487437 0.0413979668065903
3.81909547738693 0.0400224041650556
3.8391959798995 0.0383188206419581
3.85929648241206 0.0362704743081291
3.87939698492462 0.0338935590863828
3.89949748743719 0.0313129510509443
3.91959798994975 0.0289516956703636
3.93969849246231 0.0277835617679983
3.95979899497487 0.0289905278905102
3.97989949748744 0.032579621480593
4 0.0371582648419136
};
\end{axis}
\draw [gray,dashed] (2.65,0) -- (2.65,3.3);
  \draw [gray,dashed] (5.29,0) -- (5.29,3.3);
\end{tikzpicture}

%% file: figures/fl_tf_1693_new.tex
\begin{tikzpicture}

\begin{axis}[
width=\figwidth,
height=\figheight,
legend cell align={left},
log basis y={10},
tick align=outside,
tick pos=left,
x grid style={white!69.0196078431373!black},
xlabel={\(\displaystyle \omega\)},
xmin=2, xmax=16,
xtick style={color=black},
y grid style={white!69.0196078431373!black},
ylabel={$\sigma{}_{\text{max}}\text{(G(j}\omega\text{))}$},
ymin=12.8592900737875, ymax=232.786499427812,
ymode=log,
ytick style={color=black},
legend columns=3,
legend style={nodes=right,anchor=south east,at={(1.00,1.05)}}
]
\addplot [green!50.1960784313725!black]
table {%
2 14.6685447208949
2.07035175879397 16.210161289804
2.14070351758794 17.9420623547278
2.21105527638191 19.8670896886817
2.28140703517588 21.9883789456258
2.35175879396985 24.3078928243242
2.42211055276382 26.8236238960506
2.49246231155779 29.5253353004037
2.56281407035176 32.388834669999
2.63316582914573 35.3693735316936
2.7035175879397 38.3960684191085
2.77386934673367 41.371121355018
2.84422110552764 44.1788891710954
2.91457286432161 46.7079422143348
2.98492462311558 48.8820932695267
3.05527638190955 50.6869063047661
3.12562814070352 52.1756502731553
3.19597989949749 53.4495272279274
3.26633165829146 54.6242611524009
3.33668341708543 55.8040190959674
3.4070351758794 57.0779515817541
3.47738693467337 58.5403423661759
3.54773869346734 60.3209962244339
3.61809045226131 62.6040193567442
3.68844221105528 65.6165636231469
3.75879396984925 69.5881271097738
3.82914572864322 74.7049136460363
3.89949748743719 81.0877197100046
3.96984924623116 88.7980076832612
4.04020100502513 97.8511167756735
4.1105527638191 108.210557174491
4.18090452261307 119.746275295886
4.25125628140704 132.150772770127
4.321608040201 144.82604619386
4.39195979899498 156.79963778022
4.46231155778894 166.795733358202
4.53266331658291 173.590109103895
4.60301507537688 176.573748058142
4.67336683417085 176.137973851941
4.74371859296482 173.534670797673
4.81407035175879 170.343443891428
4.88442211055276 167.963399042319
4.95477386934673 167.354197378083
5.0251256281407 168.980386725697
5.09547738693467 172.842462679312
5.16582914572864 178.530590561113
5.23618090452261 185.294872605187
5.30653266331658 192.156647774531
5.37688442211055 198.084648275262
5.44723618090452 202.220023341823
5.51758793969849 204.074035861219
5.58793969849246 203.604059337325
5.65829145728643 201.136666430573
5.7286432160804 197.197883645752
5.79899497487437 192.344439071002
5.86934673366834 187.054419545567
5.93969849246231 181.684065508718
6.01005025125628 176.469144019203
6.08040201005025 171.545324805525
6.15075376884422 166.968440631338
6.22110552763819 162.724190458809
6.29145728643216 158.727052235096
6.36180904522613 154.823347411678
6.4321608040201 150.832475490018
6.50251256281407 146.647840686817
6.57286432160804 142.325249198781
6.64321608040201 138.035513814236
6.71356783919598 133.926972727922
6.78391959798995 130.058192354425
6.85427135678392 126.423441187208
6.92462311557789 122.997554757882
6.99497487437186 119.764419789799
7.06532663316583 116.738496033985
7.1356783919598 114.017863344467
7.20603015075377 112.012168820291
7.27638190954774 111.815854024551
7.34673366834171 108.822350646983
7.41708542713568 104.214969227198
7.48743718592965 101.099460510186
7.55778894472362 98.603603114229
7.62814070351759 96.4231412862965
7.69849246231156 94.4598194089772
7.76884422110553 92.6686482440482
7.8391959798995 91.0193308054038
7.90954773869347 89.4861312205267
7.97989949748744 88.0460717050459
8.05025125628141 86.6794346159769
8.12060301507538 85.3703258574771
8.19095477386935 84.1066834924938
8.26130653266332 82.8798573675172
8.33165829145729 81.6840405048904
8.40201005025126 80.5157092995844
8.47236180904522 79.3730852496416
8.5427135678392 78.2555850360195
8.61306532663317 77.1632689230536
8.68341708542714 76.0963562062081
8.75376884422111 75.0548891259774
8.82412060301507 74.0385853100299
8.89447236180905 73.0468571325265
8.96482412060302 72.0789326849329
9.03517587939698 71.1340047617125
9.10552763819095 70.2113547437042
9.17587939698493 69.310427066143
9.24623115577889 68.4308505847473
9.31658291457286 67.5724139611365
9.38693467336683 66.7350097823421
9.4572864321608 65.9185683304277
9.52763819095477 65.1230026675048
9.59798994974874 64.3481795874729
9.66834170854271 63.5939194544829
9.73869346733668 62.8600179422412
9.80904522613065 62.1462776513383
9.87939698492462 61.4525370090887
9.94974874371859 60.7786853092509
10.0201005025126 60.1246547496297
10.0904522613065 59.4903834928032
10.1608040201005 58.8757498410568
10.2311557788945 58.2804871245214
10.3015075376884 57.7040991859682
10.3718592964824 57.1458019157535
10.4422110552764 56.6045121595301
10.5125628140704 56.0788914445332
10.5829145728643 55.5674344617337
10.6532663316583 55.0685796343133
10.7236180904523 54.5808164383696
10.7939698492462 54.102770411836
10.8643216080402 53.6332569663121
10.9346733668342 53.171304273603
11.0050251256281 52.7161511826688
11.0753768844221 52.2672282827109
11.1457286432161 51.8241299328141
11.21608040201 51.3865834795161
11.286432160804 50.9544198338723
11.356783919598 50.5275476399583
11.427135678392 50.1059317555312
11.4974874371859 49.6895758077054
11.5678391959799 49.2785081465747
11.6381909547739 48.8727704631383
11.7085427135678 48.4724084967648
11.7788944723618 48.0774644849402
11.8492462311558 47.6879712037096
11.9195979899497 47.3039475587171
11.9899497487437 46.9253957124073
12.0603015075377 46.5522996988136
12.1306532663317 46.1846254115478
12.2010050251256 45.8223217938351
12.2713567839196 45.4653230238405
12.3417085427136 45.1135514759198
12.4120603015075 44.7669212361461
12.4824120603015 44.4253419314369
12.5527638190955 44.0887225700952
12.6231155778894 43.7569749776901
12.6934673366834 43.4300162707675
12.7638190954774 43.1077697298058
12.8341708542714 42.7901635499334
12.9045226130653 42.4771273851898
12.9748743718593 42.1685873345715
13.0452261306533 41.8644607640673
13.1155778894472 41.5646526679338
13.1859296482412 41.26905483698
13.2562814070352 40.9775480636272
13.3266331658291 40.6900065273193
13.3969849246231 40.4063029297245
13.4673366834171 40.126313056062
13.5376884422111 39.8499189898707
13.608040201005 39.5770108171965
13.678391959799 39.3074870726576
13.748743718593 39.0412543428045
13.8190954773869 38.7782264226564
13.8894472361809 38.5183233239557
13.9597989949749 38.2614703139692
14.0301507537688 38.0075970762882
14.1005025125628 37.7566370240813
14.1708542713568 37.5085267609492
14.2412060301508 37.2632056724372
14.3115577889447 37.0206156198522
14.3819095477387 36.7807007168059
14.4522613065327 36.5434071670788
14.5226130653266 36.3086831460215
14.5929648241206 36.0764787150408
14.6633165829146 35.846745758303
14.7336683417085 35.6194379345058
14.8040201005025 35.3945106384885
14.8743718592965 35.1719209695322
14.9447236180905 34.9516277014807
15.0150753768844 34.7335912547964
15.0854271356784 34.5177736670096
15.1557788944724 34.3041385613741
15.2261306532663 34.0926511123609
15.2964824120603 33.8832780077329
15.3668341708543 33.6759874055832
15.4371859296482 33.4707488871879
15.5075376884422 33.267533405386
15.5778894472362 33.0663132269041
15.6482412060301 32.8670618710454
15.7185929648241 32.6697540427788
15.7889447236181 32.474365562078
15.8592964824121 32.2808732891401
15.929648241206 32.0892550464692
16 31.8994895381359
};
\addlegendentry{Full model}
\addplot [blue, dotted,mark=+]
table {%
2 14.7036709433386
2.07035175879397 16.2563208757498
2.14070351758794 17.9889161749436
2.21105527638191 19.9059313963222
2.28140703517588 22.013220049337
2.35175879396985 24.3158984864538
2.42211055276382 26.8151760624487
2.49246231155779 29.5038034528873
2.56281407035176 32.359990309834
2.63316582914573 35.3403237887082
2.7035175879397 38.3736655736578
2.77386934673367 41.3600929473287
2.84422110552764 44.1803623523015
2.91457286432161 46.7192404951808
2.98492462311558 48.898155385071
3.05527638190955 50.7025121137763
3.12562814070352 52.1869778760986
3.19597989949749 53.4543572893494
3.26633165829146 54.6215542881734
3.33668341708543 55.7939417896699
3.4070351758794 57.0624801193966
3.47738693467337 58.5236187558669
3.54773869346734 60.3086164267694
3.61809045226131 62.6011475168985
3.68844221105528 65.6255365312614
3.75879396984925 69.6068311730562
3.82914572864322 74.7272474373869
3.89949748743719 81.1062060161131
3.96984924623116 88.8070157792346
4.04020100502513 97.848972763128
4.1105527638191 108.199858579789
4.18090452261307 119.732584772194
4.25125628140704 132.140137406757
4.321608040201 144.822055440999
4.39195979899498 156.801488704067
4.46231155778894 166.799326354958
4.53266331658291 173.592323844162
4.60301507537688 176.575639187781
4.67336683417085 176.142316426143
4.74371859296482 173.540745946329
4.81407035175879 170.346309321117
4.88442211055276 167.958902817951
4.95477386934673 167.344087920226
5.0251256281407 168.972066295665
5.09547738693467 172.843090324673
5.16582914572864 178.540324800237
5.23618090452261 185.305414184392
5.30653266331658 192.157560090828
5.37688442211055 198.073093472052
5.44723618090452 202.205033721126
5.51758793969849 204.070027822484
5.58793969849246 203.618313264224
5.65829145728643 201.161588717485
5.7286432160804 197.21460878348
5.79899497487437 192.335322921323
5.86934673366834 187.016570764876
5.93969849246231 181.635884694492
6.01005025125628 176.444426570161
6.08040201005025 171.573981657002
6.15075376884422 167.050482826141
6.22110552763819 162.813805849404
6.29145728643216 158.749451387272
6.36180904522613 154.732791312552
6.4321608040201 150.673512738033
6.50251256281407 146.540986814095
6.57286432160804 142.361310933771
6.64321608040201 138.194636383995
6.71356783919598 134.109546289231
6.78391959798995 130.16599932949
6.85427135678392 126.409273812131
6.92462311557789 122.87246298008
6.99497487437186 119.58563916353
7.06532663316583 116.595489033147
7.1356783919598 114.016808483753
7.20603015075377 112.190353088139
7.27638190954774 111.704386348362
7.34673366834171 108.976024389952
7.41708542713568 104.122716972166
7.48743718592965 100.960837137155
7.55778894472362 98.5540659016754
7.62814070351759 96.4669344065681
7.69849246231156 94.5653198095781
7.76884422110553 92.7979476178363
7.8391959798995 91.1395310994094
7.90954773869347 89.5744551209145
7.97989949748744 88.0912957880224
8.05025125628141 86.6807137007644
8.12060301507538 85.3345546396259
8.19095477386935 84.0454394859615
8.26130653266332 82.8065780984024
8.33165829145729 81.6117001804509
8.40201005025126 80.4550559511677
8.47236180904522 79.3314622381782
8.5427135678392 78.2363769525785
8.61306532663317 77.1659854053455
8.68341708542714 76.1172799917637
8.75376884422111 75.0881134219102
8.82412060301507 74.0772071898389
8.89447236180905 73.0841026138052
8.96482412060302 72.109051285634
9.03517587939698 71.1528532041883
9.10552763819095 70.2166611796143
9.17587939698493 69.3017762830308
9.24623115577889 68.4094594999054
9.31658291457286 67.5407796299899
9.38693467336683 66.6965088891398
9.4572864321608 65.8770683752108
9.52763819095477 65.0825179529164
9.59798994974874 64.3125804955094
9.66834170854271 63.566688901391
9.73869346733668 62.8440451388969
9.80904522613065 62.1436826908087
9.87939698492462 61.4645262176682
9.94974874371859 60.8054444481568
10.0201005025126 60.1652940416405
10.0904522613065 59.542953511847
10.1608040201005 58.937347367635
10.2311557788945 58.3474614582238
10.3015075376884 57.7723510659142
10.3718592964824 57.2111435245364
10.4422110552764 56.6630370773227
10.5125628140704 56.1272974163682
10.5829145728643 55.6032529879564
10.6532663316583 55.090289802939
10.7236180904523 54.5878462132887
10.7939698492462 54.0954079179957
10.8643216080402 53.6125033324933
10.9346733668342 53.1386993765583
11.0050251256281 52.6735976885886
11.0753768844221 52.2168312469245
11.1457286432161 51.7680613639485
11.21608040201 51.326975011825
11.286432160804 50.8932824374354
11.356783919598 50.4667150265232
11.427135678392 50.0470233818086
11.4974874371859 49.6339755855937
11.5678391959799 49.2273556231499
11.6381909547739 48.8269619483032
11.7085427135678 48.4326061767514
11.7788944723618 48.0441118957056
11.8492462311558 47.6613135805218
11.9195979899497 47.2840556103151
11.9899497487437 46.9121913753147
12.0603015075377 46.5455824691506
12.1306532663317 46.1840979595185
12.2010050251256 45.8276137308496
12.2713567839196 45.4760118928054
12.3417085427136 45.1291802486416
12.4120603015075 44.7870118177649
12.4824120603015 44.4494044071312
12.5527638190955 44.1162602264994
12.6231155778894 43.7874855429365
12.6934673366834 43.4629903703677
12.7638190954774 43.1426881903571
12.8341708542714 42.826495700682
12.9045226130653 42.5143325886246
12.9748743718593 42.2061213262386
13.0452261306533 41.9017869851567
13.1155778894472 41.6012570687835
13.1859296482412 41.3044613599705
13.2562814070352 41.0113317824934
13.3266331658291 40.7218022748523
13.3969849246231 40.4358086750902
13.4673366834171 40.1532886154796
13.5376884422111 39.874181426065
13.608040201005 39.598428046164
13.678391959799 39.3259709430368
13.748743718593 39.056754037022
13.8190954773869 38.7907226325171
13.8894472361809 38.5278233542502
13.9597989949749 38.2680040883506
14.0301507537688 38.0112139277772
14.1005025125628 37.7574031217119
14.1708542713568 37.5065230285635
14.2412060301508 37.2585260722665
14.3115577889447 37.013365701587
14.3819095477387 36.7709963521781
14.4522613065327 36.5313734111529
14.5226130653266 36.294453183961
14.5929648241206 36.0601928633795
14.6633165829146 35.828550500443
14.7336683417085 35.5994849771537
14.8040201005025 35.3729559808278
14.8743718592965 35.1489239799467
14.9447236180905 34.9273502013922
15.0150753768844 34.7081966089557
15.0854271356784 34.4914258830214
15.1557788944724 34.2770014013312
15.2261306532663 34.0648872207479
15.2964824120603 33.8550480599372
15.3668341708543 33.6474492829012
15.4371859296482 33.442056883296
15.5075376884422 33.2388374694748
15.5778894472362 33.0377582502018
15.6482412060301 32.8387870209881
15.7185929648241 32.6418921510012
15.7889447236181 32.4470425705091
15.8592964824121 32.2542077588177
15.929648241206 32.063357732668
16 31.8744630350598
};
\addlegendentry{Frequency unrestricted}
\addplot [red, dashed,mark=o]
table {%
2 16.2665701067659
2.07035175879397 17.2640836386999
2.14070351758794 18.3320794301801
2.21105527638191 19.4732245775564
2.28140703517588 20.6907126531154
2.35175879396985 21.9883216829505
2.42211055276382 23.3704809219963
2.49246231155779 24.8423512119014
2.56281407035176 26.4099230481931
2.63316582914573 28.0801366615221
2.7035175879397 29.8610292276774
2.77386934673367 31.7619156748254
2.84422110552764 33.7936114559069
2.91457286432161 35.9687081613367
2.98492462311558 38.3019160520616
3.05527638190955 40.8104915677262
3.12562814070352 43.5147725747895
3.19597989949749 46.4388492438626
3.26633165829146 49.611403006894
3.33668341708543 53.0667476575405
3.4070351758794 56.8461000216839
3.47738693467337 60.9990815370896
3.54773869346734 65.5853836178634
3.61809045226131 70.6763731750277
3.68844221105528 76.3560819624964
3.75879396984925 82.7203525045815
3.82914572864322 89.8716332624842
3.89949748743719 97.9046647597643
3.96984924623116 106.874918140733
4.04020100502513 116.738388251991
4.1105527638191 127.254131815889
4.18090452261307 137.865201435279
4.25125628140704 147.638328929472
4.321608040201 155.419198236138
4.39195979899498 160.294709316566
4.46231155778894 162.122329700931
4.53266331658291 161.647172932315
4.60301507537688 160.08810434025
4.67336683417085 158.601026403465
4.74371859296482 157.991190420571
4.81407035175879 158.690363274989
4.88442211055276 160.848762094409
4.95477386934673 164.42985450866
5.0251256281407 169.268606695771
5.09547738693467 175.093350870763
5.16582914572864 181.525584852218
5.23618090452261 188.077375082464
5.30653266331658 194.169136575118
5.37688442211055 199.1880214613
5.44723618090452 202.590994312647
5.51758793969849 204.02608732559
5.58793969849246 203.418847021457
5.65829145728643 200.977337960905
5.7286432160804 197.112513703106
5.79899497487437 192.316427619177
5.86934673366834 187.052077795014
5.93969849246231 181.686546249004
6.01005025125628 176.46997223116
6.08040201005025 171.545265578942
6.15075376884422 166.968616139272
6.22110552763819 162.72425885203
6.29145728643216 158.726826442471
6.36180904522613 154.823383056607
6.4321608040201 150.832617621668
6.50251256281407 146.647575463164
6.57286432160804 142.325148595757
6.64321608040201 138.035841307939
6.71356783919598 133.927047946135
6.78391959798995 130.057889920002
6.85427135678392 126.423359061832
6.92462311557789 122.997902963691
6.99497487437186 119.764675102681
7.06532663316583 116.738219067767
7.1356783919598 114.017529822085
7.20603015075377 112.012503551342
7.27638190954774 111.815434966155
7.34673366834171 108.822401008578
7.41708542713568 104.21483644509
7.48743718592965 101.099800498068
7.55778894472362 98.6036367616912
7.62814070351759 96.422859797571
7.69849246231156 94.4596446890455
7.76884422110553 92.6687677702962
7.8391959798995 91.0195653527364
7.90954773869347 89.4862406938002
7.97989949748744 88.0460134699619
8.05025125628141 86.6793378183307
8.12060301507538 85.3703138704266
8.19095477386935 84.1067625686059
8.26130653266332 82.8799452919829
8.33165829145729 81.6840613544771
8.40201005025126 80.5156524001244
8.47236180904522 79.3729929632745
8.5427135678392 78.2555053650457
8.61306532663317 77.1632238668345
8.68341708542714 76.0963429285481
8.75376884422111 75.054896238297
8.82412060301507 74.0386025876621
8.89447236180905 73.0468750100107
8.96482412060302 72.0789416208413
9.03517587939698 71.1340016340504
9.10552763819095 70.2113514822171
9.17587939698493 69.3104502513774
9.24623115577889 68.430926860756
9.31658291457286 67.5725491791159
9.38693467336683 66.7351774673811
9.4572864321608 65.9187183307053
9.52763819095477 65.1230875493971
9.59798994974874 64.3481842159326
9.66834170854271 63.5938752477755
9.73869346733668 62.8599879902223
9.80904522613065 62.1463084394496
9.87939698492462 61.4525829621774
9.94974874371859 60.7785218868314
10.0201005025126 60.1238038115201
10.0904522613065 59.4880798532195
10.1608040201005 58.8709773434067
10.2311557788945 58.2721026738523
10.3015075376884 57.6910431340396
10.3718592964824 57.127367680256
10.4422110552764 56.5806266546362
10.5125628140704 56.0503505455771
10.5829145728643 55.5360479594402
10.6532663316583 55.0372030611504
10.7236180904523 54.5532728341745
10.7939698492462 54.0836845956658
10.8643216080402 53.6278342593623
10.9346733668342 53.1850858419508
11.0050251256281 52.7547726341057
11.0753768844221 52.336200293021
11.1457286432161 51.9286518665185
11.21608040201 51.5313944632065
11.286432160804 51.1436869952438
11.356783919598 50.7647882063845
11.427135678392 50.39396411514
11.4974874371859 50.0304940776994
11.5678391959799 49.6736748940279
11.6381909547739 49.3228226922569
11.7085427135678 48.9772726593802
11.7788944723618 48.6363769716796
11.8492462311558 48.2995014685692
11.9195979899497 47.9660216921677
11.9899497487437 47.6353188931636
12.0603015075377 47.3067765109587
12.1306532663317 46.979777507595
12.2010050251256 46.6537028011042
12.2713567839196 46.3279309259295
12.3417085427136 46.0018389559995
12.4120603015075 45.6748046618235
12.4824120603015 45.3462098325393
12.5527638190955 45.0154446706635
12.6231155778894 44.6819131544591
12.6934673366834 44.3450392543073
12.7638190954774 44.0042738809472
12.8341708542714 43.6591024321602
12.9045226130653 43.3090527892558
12.9748743718593 42.9537035957993
13.0452261306533 42.5926926292734
13.1155778894472 42.2257250538088
13.1859296482412 41.8525813208777
13.2562814070352 41.4731244675533
13.3266331658291 41.0873065511476
13.3969849246231 40.6951739572073
13.4673366834171 40.2968713270375
13.5376884422111 39.8926438725291
13.608040201005 39.4828378808392
13.678391959799 39.0678992591273
13.748743718593 38.648370028721
13.8190954773869 38.2248827467198
13.8894472361809 37.7981529076441
13.9597989949749 37.3689694547428
14.0301507537688 36.938183605416
14.1005025125628 36.5066962636406
14.1708542713568 36.0754443501241
14.2412060301508 35.6453864246204
14.3115577889447 35.2174880018703
14.3819095477387 34.7927069713681
14.4522613065327 34.3719795213906
14.5226130653266 33.9562069402628
14.5929648241206 33.5462436247167
14.6633165829146 33.1428865692445
14.7336683417085 32.746866545137
14.8040201005025 32.3588411073099
14.8743718592965 31.9793894950348
14.9447236180905 31.609009423006
15.0150753768844 31.2481156954116
15.0854271356784 30.8970405201658
15.1557788944724 30.5560353555451
15.2261306532663 30.2252740881771
15.2964824120603 29.904857320013
15.3668341708543 29.5948175322498
15.4371859296482 29.2951248951432
15.5075376884422 29.0056935027328
15.5778894472362 28.7263878289152
15.6482412060301 28.4570292240624
15.7185929648241 28.1974022975268
15.7889447236181 27.9472610591025
15.8592964824121 27.7063347202758
15.929648241206 27.4743330825163
16 27.2509514641273
};
\addlegendentry{Frequency restricted}
\end{axis}
\draw [gray,dashed] (3,0) -- (3,3.3);
  \draw [gray,dashed] (6,0) -- (6,3.3);
\end{tikzpicture}

%% file: figures/fl_abs_err_1693_new.tex
\begin{tikzpicture}

\begin{axis}[
width=\figwidth,
height=\figheight,
legend cell align={left},
log basis y={10},
tick align=outside,
tick pos=left,
x grid style={white!69.0196078431373!black},
xlabel={\(\displaystyle \omega\)},
xmin=2, xmax=16,
xtick style={color=black},
y grid style={white!69.0196078431373!black},
ylabel={$\sigma{}_{\text{max}}\text{(G(j}\omega\text{)-}\hat{\text{G}}\text{(j}\omega\text{))}$},
ymin=0.000127950446328041, ymax=47.8480145588115,
ymode=log,
ytick style={color=black}
]
\addplot [blue, dotted,mark=+]
table {%
2 0.0794632478092895
2.07035175879397 0.0768219095553041
2.14070351758794 0.0742357606935467
2.21105527638191 0.0717903955488361
2.28140703517588 0.0695264444467953
2.35175879396985 0.0674582744361678
2.42211055276382 0.0655864419483411
2.49246231155779 0.0639051329208743
2.56281407035176 0.0624063483762042
2.63316582914573 0.0610821188662141
2.7035175879397 0.0599255187170427
2.77386934673367 0.0589309323405846
2.84422110552764 0.0580938773031431
2.91457286432161 0.057410660800622
2.98492462311558 0.0568780371644135
3.05527638190955 0.0564926965955698
3.12562814070352 0.0562499667764441
3.19597989949749 0.0561410328648119
3.26633165829146 0.0561486696955572
3.33668341708543 0.0562427740311071
3.4070351758794 0.0563781921203577
3.47738693467337 0.0564975848069557
3.54773869346734 0.0565405556398026
3.61809045226131 0.0564571275631681
3.68844221105528 0.0562206200991014
3.75879396984925 0.0558343265170721
3.82914572864322 0.0553288682725121
3.89949748743719 0.05475113843043
3.96984924623116 0.0541488453600538
4.04020100502513 0.0535557490076572
4.1105527638191 0.0529827846547039
4.18090452261307 0.0524193303837963
4.25125628140704 0.0518446472290256
4.321608040201 0.051242927613685
4.39195979899498 0.0506133402741612
4.46231155778894 0.0499716623226381
4.53266331658291 0.0493461224687871
4.60301507537688 0.0487717110070972
4.67336683417085 0.0482856363305074
4.74371859296482 0.0479243982229114
4.81407035175879 0.0477215542448708
4.88442211055276 0.0477057435677985
4.95477386934673 0.0478996628329572
5.0251256281407 0.0483206627921256
5.09547738693467 0.0489825859578423
5.16582914572864 0.0498977961255662
5.23618090452261 0.0510787144653196
5.30653266331658 0.0525390774784862
5.37688442211055 0.0542956434077768
5.44723618090452 0.05637078495944
5.51758793969849 0.0587956691309282
5.58793969849246 0.0616134174904389
5.65829145728643 0.064882399973235
5.7286432160804 0.0686808037530758
5.79899497487437 0.0731140544050601
5.86934673366834 0.0783260827889932
5.93969849246231 0.0845142243592743
6.01005025125628 0.0919446316404484
6.08040201005025 0.100957966321094
6.15075376884422 0.111940169786011
6.22110552763819 0.125208494028483
6.29145728643216 0.14074332675482
6.36180904522613 0.157752915525292
6.4321608040201 0.174332363780137
6.50251256281407 0.187870127343173
6.57286432160804 0.196462159473106
6.64321608040201 0.200095272806012
6.71356783919598 0.200280354026923
6.78391959798995 0.198835996242924
6.85427135678392 0.197127599816124
6.92462311557789 0.195973106509306
6.99497487437186 0.195850406609996
7.06532663316583 0.197196066415741
7.1356783919598 0.200867386270209
7.20603015075377 0.208868705098506
7.27638190954774 0.217036304405158
7.34673366834171 0.187107355508096
7.41708542713568 0.171416277214014
7.48743718592965 0.168132367823269
7.55778894472362 0.165808106092752
7.62814070351759 0.162754695361219
7.69849246231156 0.158747148877229
7.76884422110553 0.15381281308036
7.8391959798995 0.148082229463932
7.90954773869347 0.141757657727367
7.97989949748744 0.135072804003808
8.05025125628141 0.128247069947304
8.12060301507538 0.121454420561064
8.19095477386935 0.114814862138396
8.26130653266332 0.108403751188473
8.33165829145729 0.102269018848383
8.40201005025126 0.0964473422192482
8.47236180904522 0.0909740528721793
8.5427135678392 0.0858857359818737
8.61306532663317 0.0812175712308388
8.68341708542714 0.0769987224147606
8.75376884422111 0.0732486243625425
8.82412060301507 0.0699756585842404
8.89447236180905 0.0671782768495219
8.96482412060302 0.0648476144147448
9.03517587939698 0.0629701678890373
9.10552763819095 0.0615293490576237
9.17587939698493 0.0605057018272315
9.24623115577889 0.0598766962909788
9.31658291457286 0.0596172590537345
9.38693467336683 0.0597012432407263
9.4572864321608 0.0601029383820278
9.52763819095477 0.0607976253840426
9.59798994974874 0.0617609061492519
9.66834170854271 0.0629671684392591
9.73869346733668 0.0643876877083067
9.80904522613065 0.0659887239689134
9.87939698492462 0.0677298719179992
9.94974874371859 0.0695629662289941
10.0201005025126 0.0714319630176773
10.0904522613065 0.0732742621990542
10.1608040201005 0.075023767785566
10.2311557788945 0.076615555374431
10.3015075376884 0.0779914477622582
10.3718592964824 0.0791053346771383
10.4422110552764 0.0799269892084556
10.5125628140704 0.0804435289186016
10.5829145728643 0.0806583929499043
10.6532663316583 0.0805883935356073
10.7236180904523 0.0802597770476931
10.7939698492462 0.0797042064000258
10.8643216080402 0.0789552925448662
10.9346733668342 0.0780459621111362
11.0050251256281 0.077006696181357
11.0753768844221 0.0758645414172277
11.1457286432161 0.0746427488223601
11.21608040201 0.073360884498432
11.286432160804 0.0720352584794145
11.356783919598 0.0706795255168909
11.427135678392 0.069305328578921
11.4974874371859 0.0679228854356926
11.5678391959799 0.0665414526189802
11.6381909547739 0.0651696366863043
11.7085427135678 0.0638155517160962
11.7788944723618 0.0624868435520804
11.8492462311558 0.0611906149731514
11.9195979899497 0.0599332907278021
11.9899497487437 0.0587204614235657
12.0603015075377 0.0575567408465127
12.1306532663317 0.0564456635956474
12.2010050251256 0.0553896417838929
12.2713567839196 0.0543899908915182
12.3417085427136 0.0534470255627728
12.4120603015075 0.0525602182485848
12.4824120603015 0.0517284054924751
12.5527638190955 0.0509500181211337
12.6231155778894 0.0502233029758838
12.6934673366834 0.0495464943386526
12.7638190954774 0.0489178881746622
12.8341708542714 0.0483357812758225
12.9045226130653 0.0477982764021956
12.9748743718593 0.0473030216193778
13.0452261306533 0.0468470164113469
13.1155778894472 0.0464266165576761
13.1859296482412 0.0460377751042323
13.2562814070352 0.0456764259788517
13.3266331658291 0.0453388509332442
13.3969849246231 0.0450219095461887
13.4673366834171 0.0447231036587558
13.5376884422111 0.0444405180967209
13.608040201005 0.044172701856062
13.678391959799 0.0439185397668579
13.748743718593 0.0436771430695339
13.8190954773869 0.0434477666416769
13.8894472361809 0.0432297521034787
13.9597989949749 0.0430224910400207
14.0301507537688 0.0428254022552997
14.1005025125628 0.0426379187727822
14.1708542713568 0.0424594801539536
14.2412060301508 0.0422895289188684
14.3115577889447 0.0421275085501266
14.3819095477387 0.041972862762587
14.4522613065327 0.0418250358347068
14.5226130653266 0.0416834726893686
14.5929648241206 0.0415476196986101
14.6633165829146 0.0414169255158277
14.7336683417085 0.0412908421118543
14.8040201005025 0.0411688258378647
14.8743718592965 0.0410503390292014
14.9447236180905 0.0409348512038343
15.0150753768844 0.0408218409433313
15.0854271356784 0.0407107975161951
15.1557788944724 0.0406012228524923
15.2261306532663 0.0404926334249524
15.2964824120603 0.0403845624678926
15.3668341708543 0.0402765618260084
15.4371859296482 0.0401682040190918
15.5075376884422 0.0400590844531439
15.5778894472362 0.0399488230113793
15.6482412060301 0.0398370661167318
15.7185929648241 0.0397234881477322
15.7889447236181 0.0396077930133746
15.8592964824121 0.0394897152876108
15.929648241206 0.0393690212719921
16 0.0392455095473764
};
\addplot [red, dashed,mark=o]
table {%
2 8.39541655991422
2.07035175879397 8.78120794711511
2.14070351758794 9.20438503112801
2.21105527638191 9.66688748213683
2.28140703517588 10.1706019029944
2.35175879396985 10.7166465961035
2.42211055276382 11.3042856131364
2.49246231155779 11.9294170407529
2.56281407035176 12.5826485444094
2.63316582914573 13.2472026307994
2.7035175879397 13.8973457912166
2.77386934673367 14.498638016367
2.84422110552764 15.0116320612917
2.91457286432161 15.3998916739231
2.98492462311558 15.6407922857152
3.05527638190955 15.7346223406139
3.12562814070352 15.7069032282455
3.19597989949749 15.6022003189862
3.26633165829146 15.4724770196282
3.33668341708543 15.3652164648421
3.4070351758794 15.3154023527254
3.47738693467337 15.3432550288951
3.54773869346734 15.4577122667731
3.61809045226131 15.6638074754459
3.68844221105528 15.9706953127445
3.75879396984925 16.3972152760815
3.82914572864322 16.9734995078982
3.89949748743719 17.7385859021625
3.96984924623116 18.7337328671803
4.04020100502513 19.9890389593587
4.1105527638191 21.498703851053
4.18090452261307 23.1819384696641
4.25125628140704 24.839424674231
4.321608040201 26.1420030535801
4.39195979899498 26.702616920355
4.46231155778894 26.2320576252369
4.53266331658291 24.683680622707
4.60301507537688 22.2774263931035
4.67336683417085 19.3911378617185
4.74371859296482 16.4032590369755
4.81407035175879 13.5843557257468
4.88442211055276 11.073295269337
4.95477386934673 8.90892272789915
5.0251256281407 7.07428871341384
5.09547738693467 5.530670007898
5.16582914572864 4.23755641433635
5.23618090452261 3.16197080898872
5.30653266331658 2.28070705600867
5.37688442211055 1.5777696519469
5.44723618090452 1.03904349240994
5.51758793969849 0.647024971436106
5.58793969849246 0.378544957794335
5.65829145728643 0.206480647175828
5.7286432160804 0.103772156588828
5.79899497487437 0.0470847275186282
5.86934673366834 0.0185949381256371
5.93969849246231 0.00599015065108957
6.01005025125628 0.00165029330628446
6.08040201005025 0.00053279726054894
6.15075376884422 0.000252537603162415
6.22110552763819 0.000353586162254406
6.29145728643216 0.000362570457127579
6.36180904522613 0.000377335012595104
6.4321608040201 0.000398754349968253
6.50251256281407 0.00042960050651747
6.57286432160804 0.000463693661139823
6.64321608040201 0.000443311002872748
6.71356783919598 0.000371669002165162
6.78391959798995 0.000370348342096854
6.85427135678392 0.000445412667999372
6.92462311557789 0.000504879892336121
6.99497487437186 0.000536798209595868
7.06532663316583 0.000549157077401639
7.1356783919598 0.000526899060379501
7.20603015075377 0.000485750416897978
7.27638190954774 0.000520352825071007
7.34673366834171 0.000509227542973134
7.41708542713568 0.000454036835808541
7.48743718592965 0.000432035014396751
7.55778894472362 0.000422130459552762
7.62814070351759 0.000400107387043973
7.69849246231156 0.00038656421900347
7.76884422110553 0.000393252844295397
7.8391959798995 0.000403109878337443
7.90954773869347 0.000403235013585259
7.97989949748744 0.000397230291243715
8.05025125628141 0.00038872642841194
8.12060301507538 0.000375241462401805
8.19095477386935 0.000353077543448428
8.26130653266332 0.000321037936181884
8.33165829145729 0.000283233640756515
8.40201005025126 0.000248755942116693
8.47236180904522 0.000229272465577848
8.5427135678392 0.000229770558609945
8.61306532663317 0.000237036121548542
8.68341708542714 0.000240595343687069
8.75376884422111 0.000238891582756609
8.82412060301507 0.000235743135252181
8.89447236180905 0.000239854949456797
8.96482412060302 0.000258798181042151
9.03517587939698 0.000291285792907309
9.10552763819095 0.000330669308997473
9.17587939698493 0.000369751056254672
9.24623115577889 0.000401766390191898
9.31658291457286 0.000420576027148878
9.38693467336683 0.000420728874307673
9.4572864321608 0.000397439941238746
9.52763819095477 0.000348085634878448
9.59798994974874 0.000280785533759542
9.66834170854271 0.000247366835081133
9.73869346733668 0.000310398213417093
9.80904522613065 0.00049362904766998
9.87939698492462 0.00090678359220561
9.94974874371859 0.00171800419840091
10.0201005025126 0.00313521541599168
10.0904522613065 0.00540094804395338
10.1608040201005 0.00878788812571571
10.2311557788945 0.01359058243692
10.3015075376884 0.020113301652631
10.3718592964824 0.0286561610814284
10.4422110552764 0.0395014915697119
10.5125628140704 0.0529020061799903
10.5829145728643 0.0690717025417025
10.6532663316583 0.0881798824972378
10.7236180904523 0.110348364888791
10.7939698492462 0.135651933589534
10.8643216080402 0.164122113084709
10.9346733668342 0.195754281753323
11.0050251256281 0.23051781450629
11.0753768844221 0.268368444632183
11.1457286432161 0.309261502284346
11.21608040201 0.353164313995211
11.286432160804 0.400065992281531
11.356783919598 0.449983195586427
11.427135678392 0.502961166072558
11.4974874371859 0.55907028042843
11.5678391959799 0.618399211282279
11.6381909547739 0.681046318861653
11.7085427135678 0.747110942085467
11.7788944723618 0.816685886629189
11.8492462311558 0.889851816241042
11.9195979899497 0.96667367471143
11.9899497487437 1.04719885443465
12.0603015075377 1.13145662749783
12.1306532663317 1.21945833284028
12.2010050251256 1.31119789069512
12.2713567839196 1.40665233075019
12.3417085427136 1.50578212972332
12.4120603015075 1.60853123766369
12.4824120603015 1.71482673103819
12.5527638190955 1.82457807686983
12.6231155778894 1.93767604258512
12.6934673366834 2.05399136039439
12.7638190954774 2.17337335641983
12.8341708542714 2.29564886379277
12.9045226130653 2.4206217945013
12.9748743718593 2.54807367888452
13.0452261306533 2.67776525932785
13.1155778894472 2.80943891451817
13.1859296482412 2.94282145018798
13.2562814070352 3.07762676345414
13.3266331658291 3.21355808348457
13.3969849246231 3.35030978500662
13.4673366834171 3.48756899417664
13.5376884422111 3.62501727551416
13.608040201005 3.76233262464664
13.678391959799 3.89919186444033
13.748743718593 4.03527342840824
13.8190954773869 4.17026043147927
13.8894472361809 4.30384389679827
13.9597989949749 4.43572599628316
14.0301507537688 4.56562317438597
14.1005025125628 4.69326904386758
14.1708542713568 4.81841696041655
14.2412060301508 4.94084220870267
14.3115577889447 5.06034374862792
14.3819095477387 5.17674549432099
14.4522613065327 5.28989711616819
14.5226130653266 5.39967437102693
14.5929648241206 5.50597898500039
14.6633165829146 5.60873812254904
14.7336683417085 5.70790348739997
14.8040201005025 5.80345010719674
14.8743718592965 5.89537485930913
14.9447236180905 5.9836947944566
15.0150753768844 6.06844531732814
15.0854271356784 6.1496782771498
15.1557788944724 6.22746001943488
15.2261306532663 6.30186944240606
15.2964824120603 6.37299609719093
15.3668341708543 6.44093836191633
15.4371859296482 6.5058017155282
15.5075376884422 6.5676971296806
15.5778894472362 6.62673959022492
15.6482412060301 6.68304675713868
15.7185929648241 6.73673776448383
15.7889447236181 6.78793216032348
15.8592964824121 6.8367489823329
15.929648241206 6.88330596301643
16 6.92771885648526
};
\end{axis}
\draw [gray,dashed] (3,0) -- (3,3.3);
  \draw [gray,dashed] (6,0) -- (6,3.3);
\end{tikzpicture}

%% file: figures/fl_rel_err_1693_new.tex
\begin{tikzpicture}

\begin{axis}[
width=\figwidth,
height=\figheight,
legend cell align={left},
legend style={fill opacity=0.8, draw opacity=1, text opacity=1, draw=white!80!black},
log basis y={10},
tick align=outside,
tick pos=left,
x grid style={white!69.0196078431373!black},
xlabel={\(\displaystyle \omega\)},
xmin=2, xmax=16,
xtick style={color=black},
y grid style={white!69.0196078431373!black},
ylabel={$\frac{\sigma{}_{\text{max}}\text{(G(j}\omega\text{)-}{\hat{\text{G}}}\text{(j}\omega\text{))}}{\sigma{}_{\text{max}}\text{(G(j}\omega\text{))}}$},
ymin=7.95780852417799e-07, ymax=1.08781088012661,
ymode=log,
ytick style={color=black}
]
\addplot [blue, dotted,mark=+]
table {%
2 0.00541725503935617
2.07035175879397 0.00473912061588333
2.14070351758794 0.00413752662463495
2.21105527638191 0.00361353357103608
2.28140703517588 0.00316196317239777
2.35175879396985 0.00277515928359961
2.42211055276382 0.00244509996868834
2.49246231155779 0.00216441683966246
2.56281407035176 0.00192678585111337
2.63316582914573 0.00172697768625955
2.7035175879397 0.00156072017746535
2.77386934673367 0.00142444609694962
2.84422110552764 0.00131496917177247
2.91457286432161 0.00122914129972102
2.98492462311558 0.00116357613514623
3.05527638190955 0.00111454221048519
3.12562814070352 0.00107808846620902
3.19597989949749 0.00105035602327046
3.26633165829146 0.00102790717001926
3.33668341708543 0.00100786242536376
3.4070351758794 0.000987740284260305
3.47738693467337 0.000965105131322216
3.54773869346734 0.000937327948454871
3.61809045226131 0.000901813144639348
3.68844221105528 0.000856805309433624
3.75879396984925 0.000802354206616233
3.82914572864322 0.000740632249903521
3.89949748743719 0.000675208756963908
3.96984924623116 0.00060979797602217
4.04020100502513 0.000547318730459002
4.1105527638191 0.000489626761363668
4.18090452261307 0.0004377533268093
4.25125628140704 0.000392314370489591
4.321608040201 0.000353823976835581
4.39195979899498 0.000322789905580675
4.46231155778894 0.000299597965226853
4.53266331658291 0.000284268053770581
4.60301507537688 0.000276211563403172
4.67336683417085 0.000274135300154501
4.74371859296482 0.000276166128662481
4.81407035175879 0.000280149051555439
4.88442211055276 0.000284024637747291
4.95477386934673 0.000286217277985226
5.0251256281407 0.000285954267997763
5.09547738693467 0.00028339439972412
5.16582914572864 0.000279491576030415
5.23618090452261 0.000275661780313557
5.30653266331658 0.000273417953981657
5.37688442211055 0.00027410323758319
5.44723618090452 0.000278759660036997
5.51758793969849 0.000288109503410382
5.58793969849246 0.000302613895277794
5.65829145728643 0.00032257867809314
5.7286432160804 0.000348283675683126
5.79899497487437 0.000380120448286373
5.86934673366834 0.000418734200342765
5.93969849246231 0.000465171362841499
6.01005025125628 0.000521023843298311
6.08040201005025 0.000588520651527789
6.15075376884422 0.000670427114026727
6.22110552763819 0.0007694522472378
6.29145728643216 0.000886700312095257
6.36180904522613 0.00101892200474018
6.4321608040201 0.00115580125045202
6.50251256281407 0.00128109712671727
6.57286432160804 0.0013803746037972
6.64321608040201 0.00144959269739304
6.71356783919598 0.00149544449446939
6.78391959798995 0.00152882331088434
6.85427135678392 0.00155926462659893
6.92462311557789 0.00159330896370318
6.99497487437186 0.00163529708534251
7.06532663316583 0.00168921198332325
7.1356783919598 0.00176171856214632
7.20603015075377 0.00186469655304692
7.27638190954774 0.00194101548745945
7.34673366834171 0.00171938351263031
7.41708542713568 0.00164483354440485
7.48743718592965 0.00166303921875359
7.55778894472362 0.00168156234514746
7.62814070351759 0.00168792152163942
7.69849246231156 0.00168057857690698
7.76884422110553 0.00165981500750162
7.8391959798995 0.00162693164357061
7.90954773869347 0.00158412991816602
7.97989949748744 0.00153411505349497
8.05025125628141 0.00147955591214326
8.12060301507538 0.00142267725162287
8.19095477386935 0.00136510985061779
8.26130653266332 0.00130796256933424
8.33165829145729 0.00125200734704425
8.40201005025126 0.00119786986984596
8.47236180904522 0.00114615744853625
8.5427135678392 0.00109750295703932
8.61306532663317 0.00105254186822785
8.68341708542714 0.00101185820522217
8.75376884422111 0.000975934082583172
8.82412060301507 0.000945124198297735
8.89447236180905 0.000919660057757756
8.96482412060302 0.0008996750090377
9.03517587939698 0.000885232992293591
9.10552763819095 0.000876344706382994
9.17587939698493 0.000872966801510122
9.24623115577889 0.000874995645667523
9.31658291457286 0.000882272151591664
9.38693467336683 0.000894601550751898
9.4572864321608 0.000911775542222822
9.52763819095477 0.000933581421213853
9.59798994974874 0.000959792593126836
9.66834170854271 0.000990144482041678
9.73869346733668 0.00102430272558098
9.80904522613065 0.00106182906624163
9.87939698492462 0.00110214932066974
9.94974874371859 0.00114452897220543
10.0201005025126 0.00118806441908288
10.0904522613065 0.00123169927468897
10.1608040201005 0.00127427282010171
10.2311557788945 0.00131460046328602
10.3015075376884 0.00135157551824712
10.3718592964824 0.0013842720204322
10.4422110552764 0.00141202505169897
10.5125628140704 0.00143447074017444
10.5829145728643 0.00145154070421317
10.6532663316583 0.0014634187783807
10.7236180904523 0.00147047593431145
10.7939698492462 0.00147320009295844
10.8643216080402 0.00147213309448015
10.9346733668342 0.0014678210959343
11.0050251256281 0.00146077994037383
11.0753768844221 0.00145147435419533
11.1457286432161 0.00144030876966248
11.21608040201 0.00142762720404783
11.286432160804 0.00141371953040134
11.356783919598 0.00139883150515296
11.427135678392 0.00138317612607354
11.4974874371859 0.00136694436069546
11.5678391959799 0.00135031386139082
11.6381909547739 0.0013334549293754
11.7085427135678 0.00131653354341482
11.7788944723618 0.00129971170945701
11.8492462311558 0.00128314569541578
11.9195979899497 0.00126698285916642
11.9899497487437 0.00125135783155558
12.0603015075377 0.00123638877604106
12.1306532663317 0.00122217432950174
12.2010050251256 0.00120879168962898
12.2713567839196 0.00119629614999102
12.3417085427136 0.00118472219131986
12.4120603015075 0.00117408606170009
12.4824120603015 0.0011643895858429
12.5527638190955 0.00115562473011392
12.6231155778894 0.00114777822282026
12.6934673366834 0.00114083527000661
12.7638190954774 0.00113478123505979
12.8341708542714 0.00112960029282005
12.9045226130653 0.0011252709244849
12.9748743718593 0.00112175969387044
13.0452261306533 0.00111901635794044
13.1155778894472 0.00111697352383972
13.1859296482412 0.00111555196226542
13.2562814070352 0.00111466957241874
13.3266331658291 0.00111425027427321
13.3969849246231 0.00111422986716928
13.4673366834171 0.00111455801075647
13.5376884422111 0.00111519720047655
13.608040201005 0.00111612021585189
13.678391959799 0.00111730723680395
13.748743718593 0.0011187433345769
13.8190954773869 0.00112041654943487
13.8894472361809 0.0011223165593138
13.9597989949749 0.00112443381519276
14.0301507537688 0.00112675900476795
14.1005025125628 0.00112928274691381
14.1708542713568 0.00113199541065844
14.2412060301508 0.00113488703281771
14.3115577889447 0.00113794727193936
14.3819095477387 0.00114116539230066
14.4522613065327 0.00114453027446182
14.5226130653266 0.00114803041800584
14.5929648241206 0.00115165396342544
14.6633165829146 0.00115538871492218
14.7336683417085 0.00115922216930477
14.8040201005025 0.00116314154639271
14.8743718592965 0.00116713383567424
14.9447236180905 0.00117118583298769
15.0150753768844 0.00117528419805119
15.0854271356784 0.00117941550660043
15.1557788944724 0.00118356631459648
15.2261306532663 0.00118772322197821
15.2964824120603 0.00119187294861719
15.3668341708543 0.00119600240197651
15.4371859296482 0.00120009875352588
15.5075376884422 0.00120414952214817
15.5778894472362 0.00120814264164398
15.6482412060301 0.00121206654470769
15.7185929648241 0.00121591022986329
15.7889447236181 0.00121966333530551
15.8592964824121 0.00122331620132767
15.929648241206 0.00122685993224152
16 0.00123028644394006
};
\addplot [red, dashed,mark=o]
table {%
2 0.572341477607876
2.07035175879397 0.541710090980919
2.14070351758794 0.513005966045069
2.21105527638191 0.486577935350242
2.28140703517588 0.462544416218442
2.35175879396985 0.440871064948075
2.42211055276382 0.421430216026882
2.49246231155779 0.404040019169224
2.56281407035176 0.38848722631149
2.63316582914573 0.374538797497471
2.7035175879397 0.361947104571267
2.77386934673367 0.350453106937804
2.84422110552764 0.339791976279776
2.91457286432161 0.329706061621289
2.98492462311558 0.319969772969311
3.05527638190955 0.310427751222496
3.12562814070352 0.301038954876751
3.19597989949749 0.291905300723297
3.26633165829146 0.28325283844957
3.33668341708543 0.275342470197679
3.4070351758794 0.268324316628441
3.47738693467337 0.262097118136438
3.54773869346734 0.256257575873916
3.61809045226131 0.250204501825784
3.68844221105528 0.243394265577034
3.75879396984925 0.235632369444507
3.82914572864322 0.227207270305154
3.89949748743719 0.21875798167221
3.96984924623116 0.210970193543112
4.04020100502513 0.20428013106058
4.1105527638191 0.198674735741228
4.18090452261307 0.193592146497943
4.25125628140704 0.187962765207877
4.321608040201 0.180506226197649
4.39195979899498 0.170297695188448
4.46231155778894 0.157270555409845
4.53266331658291 0.142195201962421
4.60301507537688 0.126164996994728
4.67336683417085 0.110090614974477
4.74371859296482 0.0945243907835591
4.81407035175879 0.0797468655993888
4.88442211055276 0.0659268348489841
4.95477386934673 0.0532339365697073
5.0251256281407 0.0418645551148927
5.09547738693467 0.0319983291267926
5.16582914572864 0.0237357441154365
5.23618090452261 0.0170645348386192
5.30653266331658 0.0118689989777754
5.37688442211055 0.00796512837155558
5.44723618090452 0.00513818303073572
5.51758793969849 0.0031705403811199
5.58793969849246 0.00185922107361904
5.65829145728643 0.00102656890382092
5.7286432160804 0.000526233622137878
5.79899497487437 0.0002447938071204
5.86934673366834 9.94092423521021e-05
5.93969849246231 3.29701486716355e-05
6.01005025125628 9.35173860255638e-06
6.08040201005025 3.10586873266848e-06
6.15075376884422 1.51248704370433e-06
6.22110552763819 2.17291701533405e-06
6.29145728643216 2.28423858455182e-06
6.36180904522613 2.43719709529185e-06
6.4321608040201 2.64369028402403e-06
6.50251256281407 2.92947038637229e-06
6.57286432160804 3.25798594241136e-06
6.64321608040201 3.21157208477046e-06
6.71356783919598 2.77516167650726e-06
6.78391959798995 2.84755873807325e-06
6.85427135678392 3.52318101624685e-06
6.92462311557789 4.10479617525702e-06
6.99497487437186 4.48211756494969e-06
7.06532663316583 4.70416440213319e-06
7.1356783919598 4.62119745910035e-06
7.20603015075377 4.33658612286405e-06
7.27638190954774 4.65365872854449e-06
7.34673366834171 4.67943891990587e-06
7.41708542713568 4.35673338653203e-06
7.48743718592965 4.27336617046756e-06
7.55778894472362 4.28108554069508e-06
7.62814070351759 4.1494954603894e-06
7.69849246231156 4.09236669540713e-06
7.76884422110553 4.24364498400519e-06
7.8391959798995 4.42883808055322e-06
7.90954773869347 4.50611740708222e-06
7.97989949748744 4.51161856004702e-06
8.05025125628141 4.48464425424723e-06
8.12060301507538 4.39545543059373e-06
8.19095477386935 4.19797248907025e-06
8.26130653266332 3.87353388843677e-06
8.33165829145729 3.46742936571016e-06
8.40201005025126 3.0895330151178e-06
8.47236180904522 2.88854168710651e-06
8.5427135678392 2.93615540033578e-06
8.61306532663317 3.07187765444349e-06
8.68341708542714 3.16171963654997e-06
8.75376884422111 3.18289168818352e-06
8.82412060301507 3.18405780263126e-06
8.89447236180905 3.28357658183208e-06
8.96482412060302 3.59048298028211e-06
9.03517587939698 4.09488814643671e-06
9.10552763819095 4.70962724198287e-06
9.17587939698493 5.33471040225763e-06
9.24623115577889 5.87112956742128e-06
9.31658291457286 6.22407876964064e-06
9.38693467336683 6.30447010766749e-06
9.4572864321608 6.02925626731625e-06
9.52763819095477 5.34504891698025e-06
9.59798994974874 4.36353499911914e-06
9.66834170854271 3.88978753319623e-06
9.73869346733668 4.9379275345149e-06
9.80904522613065 7.94301873459594e-06
9.87939698492462 1.47558365584076e-05
9.94974874371859 2.82665574232061e-05
10.0201005025126 5.21452543727245e-05
10.0904522613065 9.0786909191916e-05
10.1608040201005 0.000149261591562567
10.2311557788945 0.000233192670608304
10.3015075376884 0.000348559321371781
10.3718592964824 0.000501456977078988
10.4422110552764 0.00069785057873803
10.5125628140704 0.00094334971354266
10.5829145728643 0.00124302486178786
10.6532663316583 0.00160127395844967
10.7236180904523 0.00202174265775946
10.7939698492462 0.00250730105975973
10.8643216080402 0.00306008104612771
10.9346733668342 0.00368157758075732
11.0050251256281 0.00437281192451841
11.0753768844221 0.00513454517198791
11.1457286432161 0.00596751942937158
11.21608040201 0.00687269497369856
11.286432160804 0.00785144828625024
11.356783919598 0.00890570028834273
11.427135678392 0.0100379565542564
11.4974874371859 0.0112512588674934
11.5678391959799 0.0125490651917242
11.6381909547739 0.0139350872153917
11.7085427135678 0.0154131177974235
11.7788944723618 0.016986875147815
11.8492462311558 0.018659879918981
11.9195979899497 0.0204353700822014
11.9899497487437 0.0223162498373512
12.0603015075377 0.0243050640853016
12.1306532663317 0.0264039888160568
12.2010050251256 0.0286148287420812
12.2713567839196 0.0309390154340835
12.3417085427136 0.0333776011965509
12.4120603015075 0.0359312455100199
12.4824120603015 0.0386001920634565
12.5527638190955 0.0413842354803769
12.6231155778894 0.0442826782146859
12.6934673366834 0.0472942802413114
12.7638190954774 0.0504172071541224
12.8341708542714 0.0536489854990595
12.9045226130653 0.0569864758638383
12.9748743718593 0.0604258724312423
13.0452261306533 0.0639627314064488
13.1155778894472 0.0675920219269772
13.1859296482412 0.0713081862866169
13.2562814070352 0.0751051956226226
13.3266331658291 0.078976592970733
13.3969849246231 0.0829155241159663
13.4673366834171 0.086914763120749
13.5376884422111 0.0909667413987864
13.608040201005 0.0950635873442943
13.678391959799 0.0991971798459895
13.748743718593 0.103359215689543
13.8190954773869 0.107541288403091
13.8894472361809 0.111734975081887
13.9597989949749 0.115931927337975
14.0301507537688 0.120123962723082
14.1005025125628 0.12430315339987
14.1708542713568 0.12846190923801
14.2412060301508 0.132593053108077
14.3115577889447 0.136689886537552
14.3819095477387 0.140746244455197
14.4522613065327 0.144756538217207
14.5226130653266 0.148715786505152
14.5929648241206 0.152619634207949
14.6633165829146 0.156464359704114
14.7336683417085 0.160246871326134
14.8040201005025 0.163964694030435
14.8743718592965 0.167615947517226
14.9447236180905 0.171199317112293
15.0150753768844 0.174714018853151
15.0854271356784 0.178159760141986
15.1557788944724 0.181536697337356
15.2261306532663 0.184845391507884
15.2964824120603 0.188086763498398
15.3668341708543 0.191262049256156
15.4371859296482 0.194372756267146
15.5075376884422 0.197420621771054
15.5778894472362 0.200407573252924
15.6482412060301 0.203335691622201
15.7185929648241 0.206207177307259
15.7889447236181 0.209024319423506
15.8592964824121 0.211789468057325
15.929648241206 0.214505009637916
16 0.217173345303948
};
\end{axis}
\draw [gray,dashed] (3,0) -- (3,3.3);
  \draw [gray,dashed] (6,0) -- (6,3.3);
\end{tikzpicture}

%% file: figures/compare.tex
\begin{tikzpicture}

\definecolor{color0}{rgb}{0.12156862745098,0.466666666666667,0.705882352941177}
\definecolor{color1}{rgb}{1,0.498039215686275,0.0549019607843137}

\begin{axis}[
width=\textwidth,
height=0.7\textwidth,
legend cell align={left},
legend style={fill opacity=0.8, draw opacity=1, text opacity=1, at={(0.03,0.97)}, anchor=north west, draw=white!80!black},
tick align=outside,
tick pos=left,
x grid style={white!69.0196078431373!black},
xmin=-0.535, xmax=3.535,
xtick style={color=black},
xtick={0,1,2,3},
xticklabels={BPS-606,BPS-1142,BPS-1450,BPS-1693},
y grid style={white!69.0196078431373!black},
ylabel={time(sec)},
ymin=0, ymax=8.4,
ytick style={color=black}
]
\draw[draw=none,fill=color0] (axis cs:-0.35,0) rectangle (axis cs:0,0.617524321873983);
\addlegendimage{ybar,ybar legend,draw=none,fill=color0};
\addlegendentry{Algorithm 3}

\draw[draw=none,fill=color0] (axis cs:0.65,0) rectangle (axis cs:1,1.94182591438293);
\draw[draw=none,fill=color0] (axis cs:1.65,0) rectangle (axis cs:2,3.18545939922333);
\draw[draw=none,fill=color0] (axis cs:2.65,0) rectangle (axis cs:3,7.73806366920471);
\draw[draw=none,fill=color1] (axis cs:2.77555756156289e-17,0) rectangle (axis cs:0.35,0.525713992118835);
\addlegendimage{ybar,ybar legend,draw=none,fill=color1};
\addlegendentry{Algorithm 4}

\draw[draw=none,fill=color1] (axis cs:1,0) rectangle (axis cs:1.35,0.73235023021698);
\draw[draw=none,fill=color1] (axis cs:2,0) rectangle (axis cs:2.35,0.867683792114258);
\draw[draw=none,fill=color1] (axis cs:3,0) rectangle (axis cs:3.35,1.15602835416794);
\draw (axis cs:-0.175,0.617524321873983) ++(0pt,0.5pt) node[
  scale=1,
  anchor=south,
  text=black,
  rotate=0.0
]{0.62};
\draw (axis cs:0.825,1.94182591438293) ++(0pt,0.5pt) node[
  scale=1,
  anchor=south,
  text=black,
  rotate=0.0
]{1.94};
\draw (axis cs:1.825,3.18545939922333) ++(0pt,0.5pt) node[
  scale=1,
  anchor=south,
  text=black,
  rotate=0.0
]{3.19};
\draw (axis cs:2.825,7.73806366920471) ++(0pt,0.5pt) node[
  scale=1,
  anchor=south,
  text=black,
  rotate=0.0
]{7.74};
\draw (axis cs:0.175,0.525713992118835) ++(0pt,0.5pt) node[
  scale=1,
  anchor=south,
  text=black,
  rotate=0.0
]{0.53};
\draw (axis cs:1.175,0.73235023021698) ++(0pt,0.5pt) node[
  scale=1,
  anchor=south,
  text=black,
  rotate=0.0
]{0.73};
\draw (axis cs:2.175,0.867683792114258) ++(0pt,0.5pt) node[
  scale=1,
  anchor=south,
  text=black,
  rotate=0.0
]{0.87};
\draw (axis cs:3.175,1.15602835416794) ++(0pt,0.5pt) node[
  scale=1,
  anchor=south,
  text=black,
  rotate=0.0
]{1.16};
\end{axis}

\end{tikzpicture}